\documentclass[11pt,a4paper,reqno]{amsart}
\usepackage{mathrsfs}
\usepackage{amsmath,amssymb}
\usepackage{amsfonts}
\usepackage{latexsym,bm}
\usepackage{amsthm}
\usepackage{cite}
\usepackage{amssymb}
\usepackage{amssymb,amscd}
\usepackage{ulem}
\usepackage{amsbsy}
\usepackage{fancyhdr,graphicx}
\usepackage{indentfirst}
\usepackage{graphics,color}
\usepackage{epsfig}
\usepackage{xcolor}
\usepackage{overpic}
\usepackage{caption}
\usepackage{subcaption}
\usepackage{multirow}
\usepackage{diagbox}
\usepackage{array}
\usepackage{makecell}
\usepackage{arydshln}
\usepackage{booktabs}
\usepackage{threeparttable}
\usepackage{extarrows}
\usepackage{enumerate}
\usepackage{enumitem}
\newtheorem {theorem} {Theorem}
\newtheorem {proposition} [theorem]{Proposition}

\newtheorem {lemma}  [theorem]{Lemma}

\newtheorem {definition} [theorem]{Definition}

\theoremstyle{remark}

\theoremstyle{definition}

\newtheorem {remark} [theorem]{Remark}


\newpage

\begin{document}

\title[A criterion for two limit cycles in Abel equations]
{A Chebyshev criterion for at most two non-zero limit cycles in Abel equations}

\author[J. Huang, R. Tian and Y. Zhao]
	{Jianfeng Huang$^{1}$, Renhao Tian$^{2*}$, Yulin Zhao$^3$}
\begingroup
\renewcommand\thefootnote{}
\footnotetext{$*$Corresponding author. Email: tianrh@mail2.sysu.edu.cn}
\endgroup

\address{$^1$ Department of Mathematics, Jinan University, Guangzhou 510632, P. R. China}
\email{thuangjf@jnu.edu.cn}

\address{$^2$ School  of Mathematics (Zhuhai), Sun Yat-sen University, Zhuhai Campus, Zhuhai, 519082,  China}
\email{tianrh@mail2.sysu.edu.cn}

\address{$^3$ School  of Mathematics (Zhuhai), Sun Yat-sen University, Zhuhai Campus, Zhuhai, 519082,  China}
\email{mcszyl@mail.sysu.edu.cn}

\begin{abstract}
    In this paper, we investigate the maximum number of limit cycles of the reduced Abel equation $\dot{x}=A(t)x^{3}+B(t)x^{2}$ on an interval $[0,T]$.
    The Smale-Pugh problem asks whether this maximum number is bounded in terms of a given class of coefficients.
    We establish for the first time a Chebyshev criterion, providing a positive answer to the problem when this class spanned by an extended Chebyshev system (ET-system) $\mathcal{F}=\{f_{0},f_{1},f_{2}\}$ on $[0,T)$ with $f_{0}\not=0$.
    As an application, we prove that the equation has at most three limit cycles (including $x=0$) when the coefficients $A$ and $B$ are both linear trigonometric functions or quadratic polynomials. This reestablishes the result of Yu et al. (J. Differ. Equ., 2024) and improves the work of Bravo et al. (Disc. Cont. Dyn. Syst., 2015 \& J. Differ. Equ., 2024). We also obtain the same maximum number of limit cycles for the equation with trinomial coefficients.
\end{abstract}

\date{}

\maketitle

\section{Introduction}
The study of limit cycles (i.e., isolated periodic solutions) of differential equations has been of wide interest for over a century, both for its applications in characterizing real-world periodic phenomena and for its theoretical role in understanding the dynamical behavior of differential systems.
In the qualitative theory of differential equations, one of the most important problems, mentioned in the second part of Hilbert's 16th problem, asks for the maximum number of limit cycles that planar polynomial differential systems of degree $n$ can exhibit.
Despite the considerable progress made so far, this problem remains open and is highly challenging, even for the simplest case $n=2$.
For a comprehensive overview, we refer the reader to the surveys \cite{1,2,3,4,5}.


Motivated by the Hilbert's 16th problem, a closely related line of research, initially proposed by Pugh (see e.g., \cite{Smale0,Smale,ALN}), is to investigate the {\it closed solutions} of the following Abel equations of generalized type:
\begin{equation}\label{eq00}
    \frac{dx}{dt}=\sum^{K}_{i=0}A_{i}(t)x^{i},\indent A_i\in C^{0}[0,T],\ t\in [0,T],
\end{equation}
where the {\it closed solutions} are those satisfying the boundary condition $x(0)=x(T)$. Observe that these solutions can be extended periodically if $A_i$'s are additionally $T$-periodic functions defined on $\mathbb R$. In this context, an isolated closed solution of equation \eqref{eq00} is called a limit cycle. It is known that various planar differential systems can be reduced to equations of the form \eqref{eq00} using, e.g., polar coordinates and Cherkas' transformations. For more details the reader is referred to the works on the quadratic systems (see e.g., \cite{7,ALN}), the rigid systems (see e.g., \cite{6}), the systems with homogeneous nonlinearities and some other general systems (see e.g., \cite{7,10,HL1}). In particular, a result in \cite{ALN} shows that the Hilbert's 16th problem restricted to quadratic systems follows from a complete understanding of the limit cycles of equation \eqref{eq00} with $K=3$, known as the classical Abel equation.


The fundamental progress in the study of closed solutions of equation \eqref{eq00} mainly goes back to the 1980s, when it was observed for the first time that the classical Abel equation may be the simplest yet highly non-trivial type in the family of \eqref{eq00}.
Actually, equation \eqref{eq00} is linear for $K=1$ and becomes a Riccati equation for $K=2$, where the Poincar\'{e} map of the latter is given by a M\"{o}bius transformation, see for instance Lins-Neto \cite{ALN} and Lloyd \cite{8,9}. In both cases the equation clearly has at most $K$ limit cycles.
However, for the case $K=3$, i.e., the classical Abel equation, Lins-Neto proved in \cite{ALN} that the number of limit cycles is no longer uniformly bounded without additional restrictions.
More concretely, by applying bifurcation theory, he showed that equations of the form
\begin{equation}\label{eq0}
    \frac{dx}{dt}=A(t)x^{3}+B(t)x^{2},\indent A, B\in C^{0}[0,T],\ t\in[0,T],
\end{equation}
 can exhibit at least $n$ limit cycles when the coefficients $A$ and $B$ are both trigonometric polynomials or polynomials of degree $n$. Similar results were later extended to the case $K\geq3$ in \cite{AGY,HTV,gasull2006limit}.

In view of the above facts and background, a more specific problem for equation \eqref{eq00} naturally arises (see for instance, \cite{ALN,ilyashenko2000hilbert,Smale0}): {\it Is the maximum number of limit cycles for equation \eqref{eq00} bounded in terms of a given class of coefficients $A_i$'s (especially when the class consists of trigonometric polynomials or polynomials)}? Its initial non-trivial case, i.e., the problem for the equation of the form \eqref{eq0}, is known as the Smale-Pugh problem \cite{Smale0,Smale}.

The main challenge of this problem, including the Smale-Pugh case, lies in the unknown essential mechanism by which a given class of coefficients controls the number of limit cycles of the equations. So far, there have been several lines of exploration, each from a distinct perspective.
One concerns the definite sign property of the coefficients or their linear combinations, which essentially  arises from the ``transversality''  between some curve(s) and the orbits of the equations. In this line, Pliss \cite{Pliss} (resp. Gasull and Llibre \cite{7}) proved that equation \eqref{eq00} with $K=3$ has at most three limit cycles if $A_0\equiv0$ and $A_{3}$ (resp. $A_2$) does not change sign on $[0,T]$. Later, Gasull and Guillamon \cite{gasull2006limit} extended such result to the equation with only three non-vanishing coefficients $A_{n_1},A_{n_2}$ and $A_1$, provided that either $A_{n_1}\neq0$ or $A_{n_2}\neq0$ on $[0,T]$.
For equation \eqref{eq0}, {\'A}lvarez et al. in \cite{AGG} established a significant uniqueness criterion:
{\it If there exists a linear combination of $A$ and $B$ with a constant sign, then the equation possesses at most one non-zero limit cycle.}
We remark that this is one of the key criteria applied in the recent work \cite{YHL}, in which Yu et al. solved the Smale-Pugh problem for equation \eqref{eq0} with linear trigonometric coefficients. It will also provide valuable assistance in our subsequent analysis.
For more works along this line, the reader is referred to \cite{HL1,HL2,HZ,ilyashenko2000hilbert,alvarez2015limit}.
Another line of exploration focuses on the symmetries of the coefficients, particularly the sign-changing coefficients in the trigonometric case. See for instance \cite{MBF,BFG,bravo2008nonexistence} and the references therein. It is noteworthy that the authors of \cite{BFG} provided an estimate for the number of limit cycles of equation \eqref{eq0} when the coefficients are trigonometric functions with specific symmetries. In \cite{YHL}, this estimate formed the second key result in solving the Smale-Pugh problem for the equation with linear trigonometric coefficients.

We list some additional papers where readers can find details on the other lines of exploration on Abel equations,
such as the bifurcation of limit cycles \cite{HTV,GLT,GGM20}, the center-focus problem \cite{alvarez2017centers,BFY1,BFY,BFY2,BFY3}, and the integrability \cite{CH,F,HLM,HLM1,MH}.
We also emphasize that, in addition to the motivation of Hilbert's 16th problem, the study of limit cycles of the Abel equations (especially equation \eqref{eq0}) is inherently interesting due to its key role in analyzing various oscillatory processes in biology, physics, engineering, and other fields (see, e.g., \cite{11,12,13,14,15}).


Let us continue to explore the problem.
In this paper we focus on its initial case, i.e., the Smale-Pugh problem.
Although the previous lines of exploration have yielded significant results, a theoretical unification has yet to be established.
The above-mentioned works, however, indirectly point to a relationship between the complex dynamical behavior of equation \eqref{eq0} and the sign changes of the coefficients.
Furthermore, from the bifurcation analyses in \cite{ALN,HTV}, it is observed that limit cycles of equation \eqref{eq0} can appear due to the increase in the number of sign changes of the coefficients.
Inspired by these indications, we adopt a new perspective, in contrast with the previous lines, considering how the maximum number of limit cycles of the equation is affected by the zeros of the coefficients $A$ and $B$.

To begin this line, we restrict our consideration to the case where $A$ and $B$ in equation \eqref{eq0} have at most two zeros on $[0,T)$.
A more systematic characterization of these two coefficients can be given using the theory of Chebyshev systems (see e.g., \cite{KS}).
The theory of Chebyshev systems was initially developed in approximation theory, and its pervasive influence is now observed in areas of mathematics such as numerical analysis, matrix computations and bifurcation theory.
Recall that a set of linearly independent functions $\mathcal F=\{f_{0},f_{1},\cdots,f_{m}\}$ forms an {\it extended Chebyshev system} (ET-system) on an interval $E$, if the functions in $\mathcal F$ belong to $C^{m}(E)$ and $Z(\mathcal F)=\text{dim}(\text{span}\left(\mathcal F\right))-1=m$.
Here $\text{span}\left(\mathcal F\right):=\big\{\sum_{i=0}^{m}\lambda_{i}f_{i}\big|\lambda_{i}\in \mathbb{R}\big\}$ and $Z(\mathcal F)$ represents the maximum number of zeros (counted with multiplicities) that any non-trivial function in $\text{span}\left(\mathcal F\right)$ can have on $E$.
Then, our main result is stated as follows:

\begin{theorem}\label{theorem0}
    Consider the Abel equation \eqref{eq0}.
    Assume that the coefficients $A,B\in \text{span}(\mathcal{F})$, where $\mathcal{F}=\{f_{0},f_{1},f_{2}\}$ is an ET-system on $[0,T)$ (therefore naturally each $f_i\in C^2\left([0, T)\right)$), with each $f_i$ continuous on $[0, T]$ and $f_{0}|_{[0,T]}>0$. Then the equation has at most two non-zero limit cycles, counted with multiplicities. Moreover, this upper bound is sharp.
\end{theorem}

Recall also that a set of functions $\mathcal F$ forms an ET-system with accuracy $k$ if $Z(\mathcal F)=\text{dim}(\text{span}\left(\mathcal F\right))-1+k$ (see e.g., \cite{GI}). Under this definition, the set of coefficients $\{A,B\}$ in Theorem \ref{theorem0} is clearly an ET-system with accuracy one on $[0,T)$.
We would like to stress that such a Chebyshev property of $A$ and $B$ actually yields a ``quasi-dichotomy'', which is the key to classifying the equation into types requiring independent treatment. We detail this in Subsection 2.3.
In addition, the assumption $f_{0}|_{[0,T]}>0$ in the theorem is common  {in the study of} Chebyshev systems. As will be seen in Subsection 2.2 and Section 4, this enables us to apply the theory of rotated equations when estimating the number of limit cycles.

To the best of our knowledge, Theorem \ref{theorem0} provides for the first time a criterion bounding the maximum number of limit cycles of Abel equations by the Chebyshev property of their coefficients. In the following, we illustrate several applications of this criterion, each solving the Smale-Pugh problem for equation \eqref{eq0} with a simple yet non-trivial class of coefficients.

The first application addresses equation \eqref{eq0} with linear trigonometric coefficients, that is,
\begin{equation}\label{equationa}
    \frac{dx}{dt}=(a_{0}+a_{1}\sin t+a_{2}\cos t)x^{3}+(b_{0}+b_{1}\sin t+b_{2}\cos t)x^{2},\indent t\in[0,2\pi],
\end{equation}
where $a_0, a_1, a_2, b_0, b_1, b_2\in\mathbb R$. The systematic study of equation \eqref{equationa} began with \cite{AGG,BFG}. In \cite{6}, Gasull specifically refined the Smale-Pugh problem for this class as the 6th of his proposed 33 open problems, asking {\it whether the maximum number of limit cycles for the equation is three}. This problem remained unsolved until the previously mentioned work \cite{YHL}, which proved that:
\begin{theorem}[\cite{YHL}]\label{app1}
The equation \eqref{equationa}, with parameters $a_0,a_1,a_2,b_0,b_1,b_2\in\mathbb R$, has at most three limit cycles (including $x=0$). Moreover, this upper bound is sharp.
\end{theorem}
\noindent Here, Theorem \ref{app1} can be directly obtained from Theorem \ref{theorem0} by setting $\mathcal{F}=\{f_0,f_1,f_2\}=\{1,\sin t,\cos t\}$ and $T=2\pi$. For more details see Section 5.

In parallel with the linear trigonometric class, the second application concerns the counterpart of equation \eqref{equationa}, given by equation \eqref{eq0} with quadratic polynomial coefficients:
\begin{equation}\label{main equation}
    \frac{dx}{dt}=(a_{0}+a_{1}t+a_{2}t^{2})x^{3}+(b_{0}+b_{1}t+b_{2}t^{2})x^{2}, \indent t\in[0,1],
\end{equation}
where $a_0, a_1, a_2, b_0, b_1, b_2\in\mathbb R$. Equation \eqref{main equation} was also initially studied in \cite{AGG}, under the hypothesis that some linear combination of its coefficients has a definite sign.
Bravo et al. \cite{rotated3} considered a specific case $\frac{dx}{dt}=at(t-t_{0})x^{3}+b(t-t_{1})x^{2}$ with $a,b,t_0,t_1\in\mathbb R$. They proved that the number of positive limit cycles in this case is at most two.
Recently, a subsequent study on another case $\frac{dx}{dt}=t(t-t_{0})x^{3}+(t-t_{1})(t-1)x^{2}$ was presented in \cite{rotated2}, where the authors arrived at the same conclusion as in \cite{rotated3}.
To the best of our knowledge, {these constitute the main progress currently} on equation \eqref{main equation}. Note that the set $\{1,t,t^{2}\}$ naturally forms an ET-system on $[0,1)$. Thus Theorem \ref{theorem0} enables us to obtain the next complete result for the equation, providing an answer to the corresponding Smale-Pugh problem.
\begin{theorem}\label{main theorem1}
    The equation \eqref{main equation}, with parameters $a_0,a_1,a_2,b_0,b_1,b_2\in\mathbb R$, has at most three limit cycles (including $x=0$). Moreover, this upper bound is sharp.
\end{theorem}

We remark that for a given class of coefficients, the maximum number (finite or infinite) of limit cycles of equation \eqref{eq0} is called the Hilbert number in the literature (see e.g., \cite{ilyashenko2000hilbert,AGY,HTV}).
For the class given by trigonometric polynomials (resp. polynomials) of degree at most $n$, the corresponding Hilbert number is denoted by $\mathcal H_{\mathbb T}(n)$ (resp. $\mathcal H_{\mathbb P}(n)$).
The study of $\mathcal H_{\mathbb T}(n)$ and $\mathcal H_{\mathbb P}(n)$ in terms of $n$ comes from a variant of the Smale-Pugh problem, which is more directly related to the Hilbert's 16th problem and was highlighted by Lins-Neto in \cite{ALN}.
Prior to \cite{YHL}, the only known result was $\mathcal H_{\mathbb T}(0)=\mathcal H_{\mathbb P}(0)=\mathcal H_{\mathbb P}(1)=2$, established by the uniqueness criterion in \cite{AGG}.
We will mention this again in Section 4.
In this context, our criterion simultaneously reestablishes Theorem \ref{app1} (originally presented in \cite{YHL}) and yields Theorem \ref{main theorem1}, further identifying two previously unknown Hilbert numbers
$$\mathcal H_{\mathbb T}(1)=\mathcal H_{\mathbb P}(2)=3.$$


In addition to the above applications, Theorem \ref{theorem0} also provides an opportunity to study the Hilbert number for equation \eqref{eq0} with polynomial coefficients, in terms of the number of monomials instead of the degree.
Currently, this perspective has received increasing interest in exploring the limit cycles of polynomial differential systems (particularly the planar ones). See for instance a series of recent works \cite{app3-1,app3-2,app3-3} and the 8th open problem stated in \cite{6}.
Despite such progress, little is known about equation \eqref{eq0}, except for a result in \cite{AGG} that partially shows the uniqueness of non-zero limit cycles of the equation
$\frac{dx}{dt}=(a_{0}+a_{1}t^{m_1}+a_{2}t^{m_2})x^{3}+(b_{0}+b_{1}t^{m_1}+b_{2}t^{m_2})x^{2}$.
Here, by applying the theorem, we are able to consider equation \eqref{eq0} with trinomial coefficients, and then generalize Theorem \ref{main theorem1} as below:
\begin{theorem}[Generalization of Theorem \ref{main theorem1}]\label{main theorem}
    The Abel equation
    \begin{equation}\label{main equation1}
    \frac{dx}{dt}=(a_{0}t^{m_0}+a_{1}t^{m_1}+a_{2}t^{m_2})x^{3}+(b_{0}t^{m_0}+b_{1}t^{m_1}+b_{2}t^{m_2})x^{2},\indent t\in [0,1],
\end{equation}
    with $m_0,m_1,m_2\in\mathbb{Z}_0^+$, $0\leq m_0<m_1<m_2$ and parameters $a_0,a_1,a_2,b_0,b_1,b_2\in\mathbb R$, has at most three limit cycles (including $x=0$). Furthermore, this upper bound is sharp.
\end{theorem}
Theorem \ref{main theorem} is proved in Section 5. As will be seen, its validity requires a little more analysis, but mainly depends on the fact that the ordered set $\{1,t^{\alpha_1},t^{\alpha_2}\}$ with $\alpha_1,\alpha_2\in\mathbb R^+_0$ forms an {\it extended complete Chebyshev system} (ECT-system) on $(0,1]$.



From the results presented, there appears to be a deep connection between the Chebyshev properties of the coefficients of equation \eqref{eq0} and the number of its limit cycles.
Thus, we list a particular version of the Smale-Pugh problem, which may be worthy of further investigation:
\begin{itemize}[leftmargin=0.4cm]
 \item[] {\it Let $\mathcal{F}=\{f_{0},f_{1},\ldots,f_{n}\}$ be an ET-system on $[0,T)$, with each $f_i$ continuous on $[0,T]$ and $f_{0}|_{[0,T]}>0$. Is the maximum number of limit cycles (i.e., the Hilbert number) for equation \eqref{eq0} with coefficients $A,B\in\text{span}(\mathcal{F})$ bounded in terms of $n$? }
\end{itemize}
\noindent As the cases $n=0,1,2$ are now clear in this version, the next natural case would be $n=3$.

{
We would also like to briefly show the connection between integrability and limit cycles of the Abel equation \eqref{eq0}.
The study of its integrability is well-motivated, as this integrability directly relates to that of the generalized Li\'{e}nard equations (see e.g., \cite{F,HLM,MH}), and to the center-focus problem for planar systems with homogeneous nonlinearities (see e.g., \cite{alvarez2017centers,BFY1,F}).
So far, a number of exact integrability conditions for equation \eqref{eq0} have been established in the literature.
These conditions, such as the well-known \textit{Chiellini} condition (see e.g., \cite{CH,HLM1,HLM}), require specific differential or algebraic relations of the coefficients, and thus affect the coexistence of the first integral and limit cycles. For instance, such coexistence can appear in the polynomial case, but is precluded in the trigonometric polynomial case. In the latter case, a first integral of equation \eqref{eq0} actually implies a center (see e.g., \cite{BFY1,F}).
}

{
In contrast, the coefficients of the Abel equation \eqref{eq0} considered in this work involve a more ``flexible'' class generated by an ET-system, for which the integrability conditions are typically not satisfied.
Our results show that even in this generally non-integrable class, the dynamical behavior of equation \eqref{eq0} remains strongly constrained: The number of non-zero limit cycles does not exceed two. Roughly speaking, our study (based on Chebyshev systems) and the theory of integrability provide complementary perspectives for understanding the dynamics of the Abel equations. Their potential connection would be another interesting problem for further investigation.
}

The rest of the paper is organized as follows.
Section 2 provides several auxiliary results, including the ``quasi-dichotomy'' for certain ET-systems with accuracy one.
In Section 3, we study the multiplicity of non-zero limit cycles and the Lyapunov constants of equation \eqref{eq0}.
Then, we prove Theorem \ref{theorem0} in Section 4.
The details of the proofs for Theorem \ref{app1}, Theorem \ref{main theorem1} and Theorem \ref{main theorem} are given in Section 5.
{The last part is appendix.}

\section{Preliminaries}
In this section, we introduce three tools that will play important roles later in our analysis. The first two are adaptations of some classical results and have been already established. The third is derived from the Chebyshev property of functions, which, as far as we know, may be the one first applied to the classification of Abel equations.
They are presented one by one in the following subsections.

\subsection{A formula for the second-order derivative of the Poincar\'{e} map}
Consider the one-dimensional differential equation $\frac{dx}{dt}=S(t,x)$, $t\in[0,T]$.
Let $x=x(t,x_{0})$ be the solution of the equation satisfying initial condition $x(0,x_{0})=x_{0}$.
For the Poincar\'{e} map given by $P(x_{0})=x(T,x_0)$, Lloyd \cite{8} established two well-known derivative formulas
\begin{equation}\label{YHLeqs}
    \begin{split}
        P'(x_{0})&=\text{exp}\left[\int_{0}^{T}\frac{\partial S}{\partial x}\left(t,x\left(t,x_{0}\right)\right)dt\right],\\
        P''(x_{0})&=P'(x_{0})\left[\int_{0}^{T}\frac{\partial^{2} S}{\partial x^{2}}\left(t,x\left(t,x_{0}\right)\right)\cdot\text{exp}\left(\int_{0}^{t}\frac{\partial S}{\partial x}\left(\tau,x\left(\tau,x_{0}\right)\right)d\tau\right)dt\right].
    \end{split}
\end{equation}
Based on \eqref{YHLeqs}, Yu et al. recently proposed two refined formulas for equation \eqref{eq0} with periodic coefficients (see \cite{YHL}, proof of Theorem 3.1, Eqs. (12)-(13)). In fact, the same argument also applies to the case of general coefficients, and the result can be stated as follows.
\begin{proposition}\label{YHLproposition}
    Let $x=x(t,x_{0})$ be the solution of equation \eqref{eq0} with $x(0,x_0)=x_0$. Denote by $h(t)=\int_{0}^{t}A(s)x^{2}(s,x_0)\,ds$ and $P(x_{0})=x(T,x_0)$.
    If $x=x(t,x_{0})$ is non-zero, closed and satisfies {$P'(x_{0})=1$}, then
    \begin{equation}\label{77}
        h(T)=\int_{0}^{T}A(t)x^{2}(t,x_0)\,dt=0,
    \end{equation}
and
    \begin{equation}\label{88}
        P''(x_{0})=-\frac{2}{x_{0}^{2}}\int_{0}^{T}\exp h(t)\,dx(t,x_0).
    \end{equation}
\end{proposition}
{The proof of Proposition \ref{YHLproposition} follows exactly the same argument as in the proof of Theorem 3.1 in \cite{YHL}.
For the sake of brevity and compactness of the main text, it is arranged in the Appendix.}

\subsection{Rotated differential equations}
The theory of rotated differential equations provides an effective way to trace {the evolution of} the limit cycles of one-parameter families as the parameter varies (see e.g., \cite{Han, YHL,rotated2,rotated3}).  It was initially introduced by Duff for planar vector fields \cite{rotated1}. Here we mainly recall an adaptation drawn from \cite{Han} for one-dimensional non-autonomous differential equations.

\begin{definition}[\cite{Han}]\label{definition of rotated}
    Consider a family of differential equations
    \begin{equation}\label{nah}
        \frac{dx}{dt}=S(t,x;\alpha),\indent\ t\in [0,T],\ x\in I,\ \alpha\in J,
    \end{equation}
    where $I, J$ are intervals in $\mathbb R$.
Then, we say that \eqref{nah} defines a family of rotated equations on $[0,T]\times I$ with respect to $\alpha$,
if $\frac{\partial S}{\partial \alpha}\big|_{[0,T]\times I\times J}\geq 0$ (or $\leq 0$) and does not identically vanish along any solution of \eqref{nah}.
\end{definition}
\begin{proposition}[\cite{Han}]\label{property of rotated}
    Suppose that equation \eqref{nah} defines a family of rotated equations with respect to $\alpha$.
    \begin{itemize}
    \item[(i)] If $x=\varphi(t)$ is a stable (resp. unstable) limit cycle of equation \eqref{nah}$|_{\alpha=\alpha_0}$, then there exists $\delta>0$ such that for $\alpha\in J\cap(\alpha_0-\delta,\alpha_0+\delta)$, equation \eqref{nah} has a stable (resp. unstable) limit cycle $x= x(t;\alpha)$ with $x(t;\alpha_0)=\varphi(t)$.
    Furthermore, the monotonicity of $x= x(t;\alpha)$ with respect to $\alpha$ is characterized in Table \ref{Table4}.
    \item[(ii)] If $x=\varphi(t)$ is a semi-stable limit cycle of equation \eqref{nah}$|_{\alpha=\alpha_0}$, then there exists $\delta>0$ such that for $\alpha\in J\cap(\alpha_0-\delta,\alpha_0+\delta)\backslash\{\alpha_0\}$, equation \eqref{nah} has either two limit cycles located on distinct sides of $x=\varphi(t)$, or no limit cycles located near $x=\varphi(t)$.
    Furthermore, the existence and non-existence of these limit cycles as $\alpha$ varies from $\alpha_0$ are characterized in Table \ref{Table4}.
    \end{itemize}
    \begin{table}[!ht]
        \centering
        \newcommand{\tabincell}[2]{\begin{tabular}{@{}#1@{}}#2\end{tabular}}
        \begin{threeparttable}
        \begin{tabular}{|c|c|c|c|c|}
            \hline
            \tabincell{c}{{{Behavior}}\\ of  limit cycle} &Stable&\tabincell{c}{Unstable} &\tabincell{c}{Upper-stable\\ Lower-unstable}&\tabincell{c}{Upper-unstable\\ Lower-stable}\\
            \hline
            \tabincell{c}{$\frac{\partial S}{\partial \alpha}\geq 0$}& \tabincell{c}{increasing\\ in $\alpha$} &\tabincell{c}{decreasing\\ in $\alpha$} &\tabincell{c}{splits as $\alpha>\alpha_0$;\\ disappears as $\alpha<\alpha_0$}&\tabincell{c}{disappears as $\alpha>\alpha_0$;\\ splits as $\alpha<\alpha_0$}\\
           \hline
           $\frac{\partial S}{\partial \alpha}\leq 0$ & \tabincell{c}{decreasing\\ in $\alpha$} &\tabincell{c}{increasing\\ in $\alpha$} &\tabincell{c}{disappears as $\alpha>\alpha_0$;\\ splits as $\alpha<\alpha_0$}&\tabincell{c}{splits as $\alpha>\alpha_0$;\\ disappears as $\alpha<\alpha_0$}\\
           \hline
        \end{tabular}
        \end{threeparttable}
        \caption{Behavior of limit cycle(s) of \eqref{nah} as $\alpha$ varies. Here ``split'' means that $x=\varphi(t)$ bifurcates into two limit cycles  located on distinct sides of it, and ``disappear'' means that there are no limit cycles located near $x=\varphi(t)$.}\label{Table4}\vspace{-0.5cm}
    \end{table}
    \end{proposition}
   We remark that, under the assumption of Theorem \ref{theorem0}, equation \eqref{eq0} can be rewritten as
\begin{equation*}
    \frac{dx}{dt}=S(t,x;\lambda_{0},\mu_{0})
    =\left(\lambda_{0}f_{0}+\lambda_{1}f_{1}+\lambda_{2}f_{2}\right)x^{3} +  \left(\mu_{0}f_{0}+\mu_{1}f_{1}+\mu_{2}f_{2}\right)x^{2},
\end{equation*}
where $\lambda_{i}$ and $\mu_{i}$ are real numbers, $i=0,1,2$. Then one has
\begin{equation*}
    \text{sgn}\left(\frac{\partial S}{\partial \lambda_{0}}\right)
    =\text{sgn}(f_0)\cdot\text{sgn}(x)\neq0\ \text{and}\
    \text{sgn}\left(\frac{\partial S}{\partial \mu_{0}}\right)
    =\text{sgn}\left(f_{0}\right)>0, \ \text{ for } x\not=0.
\end{equation*}
According to Definition \ref{definition of rotated}, the equation with $A,B\in \text{span}(\mathcal{F})$, defines a family of rotated equations on each connected component of $[0,T]\times \mathbb{R}\backslash\{0\}$ with respect to both $\lambda_{0}$ and $\mu_0$.

\subsection{A quasi-dichotomy for ET-systems with accuracy one}
We begin this subsection by providing two characterizations for a set of functions  $\{A,B\}$ on an interval $E$:
\begin{itemize}
    \item[$\bm{(D.1)}$] There exists $(\lambda, \mu)\in\mathbb{R}^{2}\setminus\{(0,0)\}$ such that $\lambda A(t)+\mu B(t) \geq 0$ (or $\leq 0$) on $E$;
    \item[$\bm{(D.2)}$] $A(t)B'(t)-A'(t)B(t)>0$ (or $<0$) on $E$.
\end{itemize}
{The} following result presents a quasi-dichotomy for  {such a} set when it forms an ET-system with accuracy one.

\begin{proposition}\label{lemma1}
    Suppose that a set of non-vanishing smooth functions $\{A,B\}$ forms an ET-system with accuracy one on the interval $E$.
    Then, one of characterizations $\bm{(D.1)}$ and $\bm{(D.2)}$ must hold.
\end{proposition}
\begin{proof}
    It is sufficient to prove that $\bm{(D.1)}$ must hold when $\bm{(D.2)}$ is invalid. In this situation, there exists $t_{0}\in E$ such that $A(t_{0})B'(t_{0})-A'(t_{0})B(t_{0})=0$.
    The subsequent argument is divided into three cases:

    \noindent{Case 1}: $\left(A'(t_{0}),B'(t_{0})\right)=(0,0)$ and $A(t_{0})B(t_{0})=0$.
    We know that either $A(t)$ or $B(t)$ has a zero at $t_0$ with multiplicity at least two.
    Since $\{A,B\}$ forms an ET-system with accuracy one on $E$, this zero is unique and necessarily of multiplicity two.
    Thus,  {either} $A(t)$ {or} $B(t)$ does not change sign on the interval, which implies that $\bm{(D.1)}$ holds for either $(\lambda,\mu)=(1,0)$ or $(\lambda,\mu)=(0,1)$.

    \noindent{Case 2}: $\left(A'(t_{0}),B'(t_{0})\right)=(0,0)$ and $A(t_{0})B(t_{0})\not=0$.
    By taking $\lambda=B(t_{0})$ and $\mu=-A(t_{0})$,
    we have
    $\lambda A(t_{0})+\mu B(t_{0})=0$ and $\lambda A'(t_{0})+\mu B'(t_{0})=0$.
    Similar to case 1, it follows from assumption that $\lambda A(t)+\mu B(t)$ has a unique zero at $t=t_{0}$, necessarily of multiplicity two.
    Consequently, $\lambda A(t)+\mu B(t) \geq 0$ (or $\leq 0$) on $E$.
    Characterization $\bm{(D.1)}$ is satisfied.

    \noindent{Case 3}: $\left(A'(t_{0}),B'(t_{0})\right)\not=(0,0)$.
    Set $\lambda=B'(t_{0})$ and $\mu=-A'(t_{0})$.
    It is easy to see that $t=t_{0}$ is a double zero of $\lambda A(t)+\mu B(t)$.
    Then, the same argument as in the above two cases ensures that $\lambda A(t)+\mu B(t) \geq 0$ (or $\leq 0$) on $E$.
    Characterization $\bm{(D.1)}$ holds.

    The proof is finished.
\end{proof}

Proposition \ref{lemma1} will be utilized in the next Subsection 3.1 to classify the equation in Theorem \ref{theorem0} into two types, which exhibit different dynamical behaviors. The remark below provides additional details on the proposition to facilitate {a better} understanding.

\begin{figure}[t]
    \centering
    \begin{subfigure}[b]{0.7\textwidth}
        \includegraphics[scale=1.1]{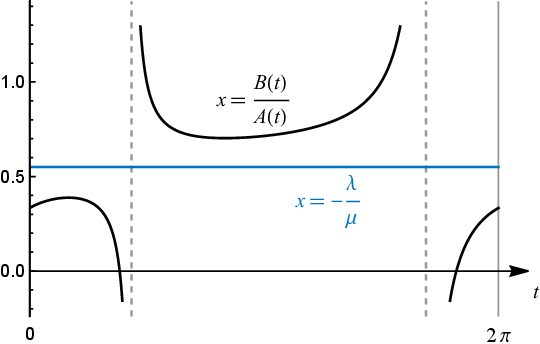}
        \caption{{Geometric interpretation of $\bm{(D.1)}$.
        Here, we take $E=(0,2\pi)$, $A(t)=1-0.5\sin\,t-2.5\cos\,t$,
        $B(t)=1-0.5\sin\,t-1.5\cos\,t$, $\mu=-1$ and $\lambda=0.55$, which implies $\lambda A+\mu B=-0.45 + 0.125 \cos t + 0.225 \sin t<0$.}}
    \end{subfigure}
    \begin{subfigure}[b]{0.7\textwidth}
            \includegraphics[scale=1.1]{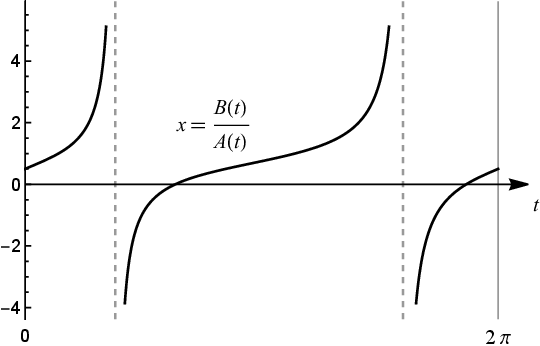}
            \caption{{Geometric interpretation of $\bm{(D.2)}$.
            Here, we take $E=(0,2\pi)$, $A(t)=1+0.1\sin\,t-3\cos\,t$ and
        $B(t)=1-2\sin\,t-2\cos\,t$, which implies $AB'-A'B=6.2 - 2.1 \cos t - \sin t>0$.}}
    \end{subfigure}
    \caption{
        Geometric interpretations of $\bm{(D.1)}$ and $\bm{(D.2)}$.
    {The black curve represents the graph for $x={B(t)}/{A(t)}$,
    the blue curve represents the graph for $x=-{\lambda}/{\mu}$,
    and the dashed lines denote the vertical lines passing through the zeros of $A(t)$.
    In subfigure $(A)$, the horizontal line $x=-{\lambda}/{\mu}$ does not cross the graph of the function $x={B(t)}/{A(t)}$.
    In subfigure $(B)$, the function $x={B(t)}/{A(t)}$ is strictly increasing on each connected components of $\{t\in (0,2\pi)|A(t)\not =0\}$.
    (For interpretation of the colors in the figures, the reader is referred to the web version of this article.)
    }}\label{figure0}
\end{figure}

\begin{remark}\label{remark1}
    We have the following comments for Proposition \ref{lemma1}.
    \begin{itemize}
        \item[(i)]
        We refer to the proposition as a ``quasi-dichotomy'' rather than a ``dichotomy'' because, under the assumption, characterizations $\bm{(D.1)}$ and $\bm{(D.2)}$ can hold simultaneously.
        To illustrate, consider $A(t)=\frac{1}{2}t^{2}+\frac{1}{4}t$, $B(t)=t^{2}-1$ and $E=[0,1]$.
        Then the set $\{A,B\}$ forms an ET-system with accuracy one on $E$.
        Taking $\lambda=1$ and $\mu=-\frac{1}{2}$, one has
        \begin{equation*}
            \lambda A(t)+\mu B(t)=\frac{1}{4}t+\frac{1}{2}\geq \frac{1}{2}, \text{ for } t\in [0,1],
        \end{equation*}
        and
        \begin{equation*}
            A(t)B'(t)-A'(t)B(t)=t^{2}+4t+1>0, \text{ for } t \in [0,1].
        \end{equation*}
        Thus, both $\bm{(D.1)}$ and $\bm{(D.2)}$ are satisfied.
        \item[(ii)]
        Both $\bm{(D.1)}$ and $\bm{(D.2)}$ have clear geometric interpretations.
        In fact, if $\bm{(D.1)}$ holds, then the horizontal line $x=-\frac{\lambda}{\mu}$ (if it exists) does not cross the graph of the function $x=\frac{B(t)}{A(t)}$ in any connected component of $\{t\in E|A(t)\not =0\}$, as illustrated in Fig. \ref{figure0}(A).
        If $\bm{(D.2)}$ is satisfied, then
        \begin{equation}\label{remeq}
            \left(\frac{B(t)}{A(t)}\right)'=\frac{A(t)B'(t)-A'(t)B(t)}{A^{2}(t)}>0\ (\text{or } <0) \text{ on } \{t\in E|A(t)\not =0\}.
        \end{equation}
        This indicates that the function $x=\frac{B(t)}{A(t)}$ has the same strict monotonicity in all connected components of $\{t\in E|A(t)\not =0\}$, as illustrated in Fig. \ref{figure0}(B).
        Note that $\bm{(D.2)}$ also ensures that $A(t)$ and $B(t)$ have no common zeros in $E$.
        Thus,
        \begin{equation}\label{remeq'}
            \lim_{t\rightarrow t_{0}}\frac{B(t)}{A(t)}=\infty\ \text{ for }\ A(t_0)=0.
        \end{equation}
    \end{itemize}
\end{remark}

\section{Multiplicity and stability of limit cycles of equation \eqref{eq0}}
The aim of this section is to study the multiplicity and stability of limit cycles of equation \eqref{eq0} under the assumption of Theorem \ref{theorem0},
including both non-zero limit cycles and {the trivial solution} $x=0$.
We recall that the set of the coefficients $\{A,B\}$ from the theorem forms an ET-system with accuracy one on $[0,T)$, and therefore Proposition \ref{lemma1} is applicable.

\subsection{Maximum multiplicity of non-zero limit cycles}
For the sake of brevity, in the following we will use the notation $x=x(t)$ to represent a non-zero limit cycle of equation \eqref{eq0}, with initial value $x(0)=x_0$. Also we define $\varphi(t):=-\frac{B(t)}{A(t)}$.
Note that equation \eqref{eq0} can be rewritten as
$\frac{dx}{dt}=A(t)x^{2}\left(x-\varphi(t)\right)$ for $t\in E\setminus\{t|A(t)=0\}$.
Thus, $x=\varphi(t)$ is the horizontal isocline of the equation.

We shall begin with the following auxiliary result.
\begin{lemma}\label{lemma2}
    For equation \eqref{eq0}, assume that the set of coefficients $\{A,B\}$ forms an ET-system with accuracy one on $[0,T)$, and satisfies $\bm{(D.2)}$ with $E=[0,T)$. If $x=x(t)$ is a non-zero limit cycle of the equation, then the following statements hold.
    \begin{itemize}
        \item[(i)] $x(t)-\varphi(t)$ has at most one zero in each connected component of $\{t\in(0,T)|A(t)\not= 0\}$.
        Furthermore, all these zeros are simple  if they exist, and their total number is either one or two.
        \item[(ii)] Each zero of $x(t)-\varphi(t)$ is a local extremum point of $x(t)$, and $x(t)$ has no other stationary point in $(0,T)$.
    \end{itemize}
\end{lemma}
\begin{proof}
    Denote by $V=\{t\in[0,T)|A(t)\not= 0\}$. By assumption and inequality \eqref{remeq} in Remark \ref{remark1}, we know that $\varphi'(t)$ has definite sign on $V$. In the following we prove the case $\varphi'(t)>0$, and the opposite case follows from a similar argument.

    (i) First, for each zero $t_{0}\in V$ of $x(t_0)-\varphi(t_0)$, we have
    \begin{equation}\label{xvarph}
        \begin{split}
            \big(x(t)-\varphi(t)\big)'\big|_{t=t_{0}}
                &=A(t_{0})x_{0}^{2}(t_{0})\left(x(t_{0})-\varphi(t_{0})\right)-\varphi'(t_{0})\\
                &=-\varphi'(t_{0})
                <0.
        \end{split}
    \end{equation}
Hence, $x(t)-\varphi(t)$ has at most one zero in each connected component of $V$, which must be simple if it exists.

Since $\{A,B\}$ is an ET-system with accuracy one, $A(t)$ has at most two zeros in $(0,T)$. From the above argument, the number of positive zeros of $x(t)-\varphi(t)$ in $V$ is at most three. Next, we show by contradiction that this upper bound cannot be achieved. Actually, if $x(t)-\varphi(t)$ has exactly three positive zeros in $V$, namely $t_1<t_2<t_3$, then it follows from \eqref{xvarph} that $x(t)-\varphi(t)>0$ (resp. $<0$) on $[0,t_1)$ (resp. $(t_3,T)$). We get
\begin{align}\label{lemma2-eq1}
\varphi(0)< x(0)=x(T)\leq\varphi(T).
\end{align}
Furthermore, in this case $A(t)$ has two zeros $a$ and $b$, satisfying $t_1<a<t_2<b<t_3$.
Note that $\varphi(t)$ is monotonically increasing on $[0,a)$, $(a,b)$ and $(b,T)$. This together with \eqref{remeq'} in Remark \ref{remeq} implies that
 \begin{align}\label{lemma2-eq2}
   \lim_{t\rightarrow a^{+}}\varphi(t)=\lim_{t\rightarrow b^{+}}\varphi(t)=-\infty\ \text{ and }\
   \lim_{t\rightarrow b^{-}}\varphi(t)=\lim_{t\rightarrow a^{-}}\varphi(t)=+\infty.
 \end{align}
Thus, according to \eqref{lemma2-eq1} and \eqref{lemma2-eq2}, there exists $\lambda\in (\varphi(0),\varphi(T))$ such that $\varphi(t)-\lambda=-\frac{B(t)+\lambda A(t)}{A(t)}$ has a zero in each of $(0,a)$, $(a,b)$ and $(b,T)$. This contradicts the assumption that $\{A, B\}$ forms an ET-system with accuracy one on $[0,T)$.
Consequently, $x(t)-\varphi(t)$ has at most two positive zeros in $V$.

Finally, since $x(t)$ is a non-zero limit cycle, there exists $t_{0}\in (0,T)$ such that $x'(t_0)=0$, which yields
\begin{align}\label{lemma2-eq3}
  A(t_{0})x(t_{0})+B(t_{0})=\frac{x'(t_{0})}{x^{2}(t_{0})}=0.
\end{align}
 As characterization $\bm{(D.2)}$ ensures that $A(t)$ and $B(t)$ have no common zeros in $[0,T)$, we get $A(t_{0})\neq0$ and therefore $x(t_0)-\varphi(t_0)=0$. Hence, $x(t)-\varphi(t)$ has at least one zero in $\{t\in[0,T)|A(t)\not= 0\}$. This completes the proof of statement (i).

(ii)
 It is clear from \eqref{lemma2-eq3} and the following analysis that any stationary point of $x(t)$ in $(0,T)$ must be the zero of $x(t)-\varphi(t)$ in $V$. Conversely, for each positive zero $t_{0}\in V$ of $x(t)-\varphi(t)$, we have $x'(t_0)=0$ and
\begin{equation}
    \begin{split}
        x''(t_{0})&=\left[A(t)x^{2}(t)\left(x(t)-\varphi(t)\right)\right]'\big|_{t=t_0}\\
        &=A(t_{0})x^{2}(t_{0})\left(x'(t_{0})-\varphi'(t_{0})\right)\\
        &=-A(t_{0})x^{2}(t_{0})\varphi'(t_{0})\not=0.
    \end{split}
\end{equation}
Hence, $t=t_{0}$ is a local extremum point of $x(t)$.
\end{proof}

By Proposition \ref{lemma1}, Lemma \ref{lemma2}, and the uniqueness criterion mentioned in Section 1 and established in \cite[Theorem A] {AGG}, we are able to give the following characterization for the multiplicity and stability of non-zero limit cycles of equation \eqref{eq0}.
\begin{proposition}\label{multiplicity of non-zero limit cycle}
    For equation \eqref{eq0}, assume that the set of coefficients $\{A,B\}$ forms an ET-system with accuracy one on $[0,T)$.
    Then, one of the following statements holds.
    \begin{itemize}
        \item[(i)] The set $\{A,B\}$ satisfies $\bm{(D.1)}$ with $E=[0,T)$, and equation \eqref{eq0} has at most one non-zero limit cycle.
        Moreover, if this limit cycle exists, it is hyperbolic.
        \item[(ii)] The set $\{A,B\}$ satisfies $\bm{(D.2)}$ with $E=[0,T)$, and the multiplicity of each non-zero limit cycle of equation \eqref{eq0} is at most two.
        Moreover, if a non-hyperbolic limit cycle exists, then it must be lower-stable and upper-unstable (resp. lower-unstable and upper-stable) when $A'(t)B(t)-A(t)B'(t)>0$ (resp. $<0$).
    \end{itemize}
\end{proposition}
Before starting the proof of Proposition \ref{multiplicity of non-zero limit cycle},
we stress that statement (i) with $A(t)$ and $B(t)$ being linearly independent is already established in \cite[Theorem A] {AGG}.
Also, the statement is trivial when $A(t)$ and $B(t)$  {are} linearly dependent. In fact, if $\lambda A(t)+\mu B(t) \equiv 0$ for certain $(\lambda,\mu)\neq(0,0)$, then equation \eqref{eq0} can be written as
\begin{align*}
  \frac{dx}{dt}=B(t)\left(-\frac{\mu}{\lambda}x^{3}+x^{2}\right) \ \
  \left(\text{resp. } \frac{dx}{dt}=A(t)\left(x^{3}-\frac{\lambda}{\mu}x^{2}\right)\right)
\end{align*}
when $\lambda\not=0$ (resp. $\mu\not=0$). Therefore, the equation is integrable, which only yields either a center or a unique and hyperbolic non-zero limit cycle.

The proof of statement (ii) is quite long, but essentially consists of four analogous parts, each dealing with an independent case. Here we provide an outline.
In the first step, we divide the argument for the non-zero limit cycle $x=x(t)$ into two cases according to Lemma \ref{lemma2}:
\begin{itemize}
  \item $x(t)-\varphi(t)$ has one zero in $(0,T)$ (i.e. $x(t)$ has one extremum point in $(0,T)$);
  \item $x(t)-\varphi(t)$ has two zeros in $(0,T)$ (i.e. $x(t)$ has two extremum points in $(0,T)$).
\end{itemize}
Then, due to the assumption that $A(t)$ has either one or two simple zero(s) in
$(0,T)$, each of the above cases is further subdivided into two subcases. All these four subcases will be verified one by one, using the formula in Proposition \ref{YHLproposition} and the approach developed from \cite[Theorem 3.1] {YHL}.


Now, we present the proof of Proposition \ref{multiplicity of non-zero limit cycle}.
\begin{proof}[Proof of Proposition \ref{multiplicity of non-zero limit cycle}]
First, from assumption and Proposition \ref{lemma1}, the set of coefficients $\{A,B\}$ satisfies either $\bm{(D.1)}$ or $\bm{(D.2)}$, with $E=[0,T)$. Therefore it is sufficient to verify the assertion on the limit cycles of equation \eqref{eq0} in each statement.

    (i) By the above discussion, the validity of statement (i) is already known and here we omit the proof.

    (ii) Taking change of variable $x\mapsto -x$ if necessary, we may restrict our consideration to the case $A'(t)B(t)-A(t)B'(t)>0$, i.e., the function $\varphi(t)=-\frac{B(t)}{A(t)}$ is strictly increasing in each connected component of $\{t\in[0,T)|A(t)\not=0\}$.

    We also exclude the cases where $A(t)$ has either no zeros or a double zero in $(0,T)$, because both are reduced to those in statement (i).
    Let $Z(A):=\{t\in(0,T)|A(t)=0, A'(t)\not=0\}$. Then, by assumption the set $Z(A)$ contains either one or two elements.

    Now assume that $x=x(t)$ is a non-zero limit cycle of equation \eqref{eq0} with initial value $x(0)=x_0$.
    As stated in the outline, and based on Lemma \ref{lemma2}, we divide the argument into two cases:
    $x(t)-\varphi(t)$ has exactly one zero in $(0,T)\backslash Z(A)$, and $x(t)-\varphi(t)$ has two zeros in $(0,T)\backslash Z(A)$.
    \vskip0.2cm

\noindent\textbf{ \textbf{Case 1}}: $x(t)-\varphi(t)$ has exactly one zero in $(0,T)\backslash Z(A)$ (i.e., $x(t)$ has exactly one extremum point), denoted by $t=t_{1}$.
    Clearly, $x(t)$ is strictly monotonic on $[0,t_1]$ and $[t_1,T]$, respectively.
    Let $h(t)=\int_{0}^{t}A(s)x^{2}(s)\,ds$.
    From Proposition \ref{YHLproposition}, if $x=x(t)$ is non-hyperbolic, then its multiplicity and stability can be determined by
\begin{equation}\label{express}
    \begin{split}
        P''(x_{0})&=-\frac{2}{x_{0}^{2}}\int_{0}^{T}\exp h(t)\,dx(t)\\
        &=-\frac{2}{x_{0}^{2}}\left(\int_{0}^{t_{1}}\exp h(t)\,dx(t)+\int_{t_1}^{T}\exp h(t)\,dx(t)\right)\\
        &=-\frac{2}{x_{0}^{2}}\left(\int_{x_{0}}^{x_{1}}\exp h(\tau_{1}(x))\,dx+\int_{x_{1}}^{x_{0}}\exp h(\tau_{2}(x))\,dx\right)\\
        &=-\frac{2}{x_{0}^{2}}\left(\int_{x_0}^{x_1}\left(\exp h(\tau_{1}(x))-\exp h(\tau_{2}(x))\right)\,dx\right),
    \end{split}
\end{equation}
where $x_{1}=x(t_{1})$, and $\tau_{1}$ and $\tau_{2}$ represent the inverse functions of $x|_{[0,t_{1}]}$ and $x|_{[t_{1},T]}$, respectively.
For convenience, we assume that $x_{0}<x_{1}$. The case $x_{0}>x_{1}$ can be handled analogously and is explained later at the end of Case 1.
Thus, we obtain that
\begin{equation}\label{eq11-1}
    \tau_{1}'|_{(x_{0},x_{1})}>0,\ \ \tau_{2}'|_{(x_{0},x_{1})}<0,\ \ \tau_{1}(x_{1})=\tau_{2}(x_{1})=t_{1},\ \ \tau_{1}(x_{0})=0,\ \ \tau_{2}(x_{0})=T.
\end{equation}

Next, we determine the sign of $P''(x_{0})$ by \eqref{express}. This can be done once the sign of the function
\begin{equation*}
    W(s):=h(\tau_{1}(s))-h(\tau_{2}(s))=\int_{\tau_{2}(s)}^{\tau_{1}(s)}A(t)x^2(t)dt
\end{equation*}
is clarified on $(x_{0},x_{1})$.
Note that we have $W(x_{0})=W(x_{1})=0$ by Proposition \ref{YHLproposition} and \eqref{eq11-1}.
Thus, the problem is reduced to analyzing the first-order derivative $W'(s)$.
Clearly, $W'(s)$ is continuous on $(x_{0},x_{1})$, and we get
\begin{equation*}
    \begin{split}
        W'(s)&=A\left(\tau_{1}(s)\right)x^2\left(\tau_{1}(s)\right)\tau_{1}'(s)-A\left(\tau_{2}(s)\right)x^2\left(\tau_{2}(s)\right)\tau_{2}'(s)\\
        &=\frac{A\left(\tau_{1}(s)\right)}{A\left(\tau_{1}(s)\right)x\left(\tau_{1}(s)\right)+B\left((\tau_{1}(s))\right)}-\frac{A\left(\tau_{2}(s)\right)}{A\left(\tau_{2}(s)\right)x\left(\tau_{2}(s)\right)+B\left((\tau_{2}(s))\right)}.
    \end{split}
\end{equation*}
Recall that $Z(A)$ is the set of the simple zeros of $A(t)$ in $(0,T)$.
For $s\in (x_{0},x_{1})\setminus\{x(t)|t\in Z(A)\}$, we obtain
\begin{equation}\label{eq011}
    W'(s)=\frac{1}{x(\tau_{1}(s))-\varphi(\tau_{1}(s))}-\frac{1}{x(\tau_{2}(s))-\varphi(\tau_{2}(s))}.
\end{equation}
Moreover, since $x\circ\tau_{1}=x\circ\tau_{2}=\text{id}$, we can write \eqref{eq011} as
\begin{equation}\label{case1W'}
    W'(s)=\frac{\varphi(\tau_{1}(s))-\varphi(\tau_{2}(s))}{\left(x(\tau_{1}(s))-\varphi(\tau_{1}(s))\right)\left(x(\tau_{2}(s))-\varphi(\tau_{2}(s))\right)}.
\end{equation}
The following discussion is subdivided into two subcases:

\begin{figure}[h]
    \centering
    \begin{subfigure}[b]{0.8\textwidth}
        \includegraphics[scale=0.65]{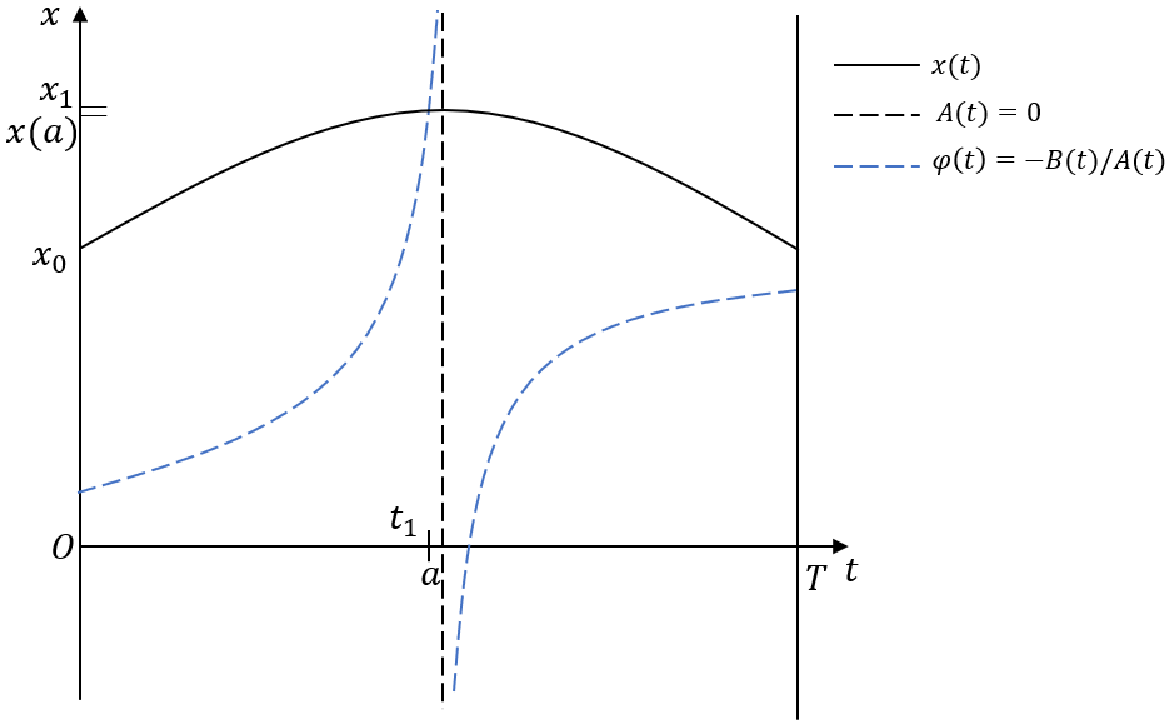}
        \caption{ {Subcase where $Z(A)$ is a one-element set $\{a\}$.}}
    \end{subfigure}
    \begin{subfigure}[b]{0.8\textwidth}
            \includegraphics[scale=0.65]{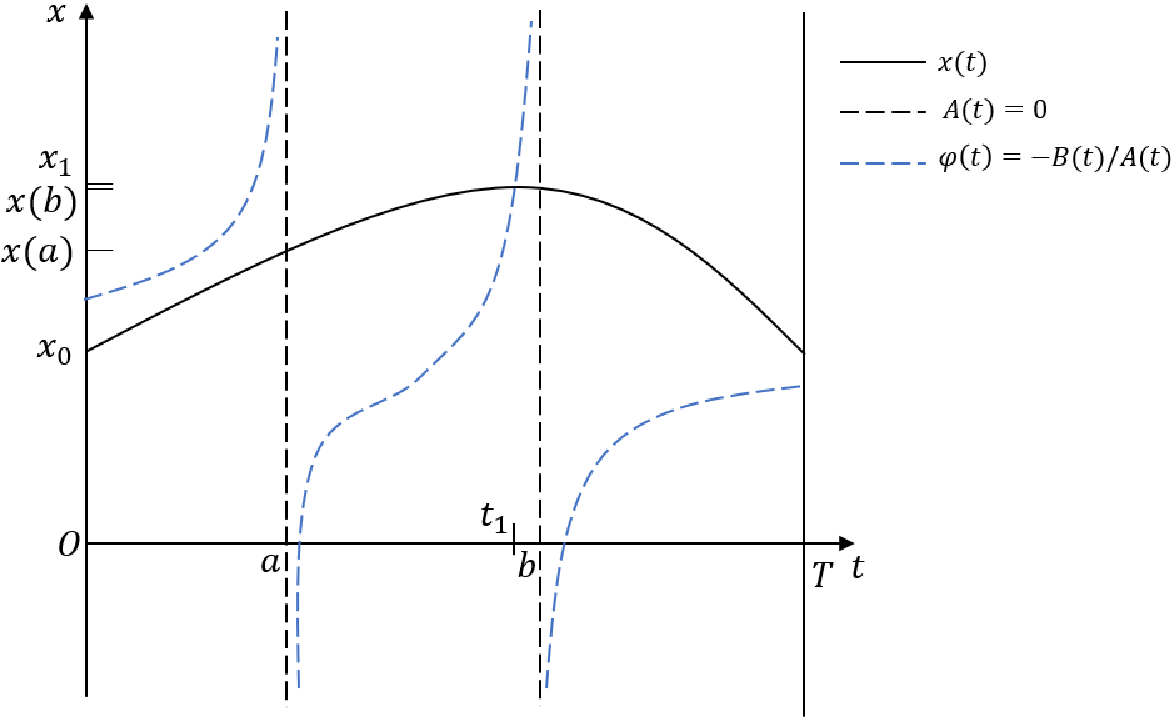}
            \caption{{Subcase where $Z(A)$ is a two-element set $\{a,b\}$ with $a<b$.}}
    \end{subfigure}
    \caption{Case where $x(t)-\varphi(t)$ has exactly one zero $t=t_{1}$ in $[0,T)$.
    {The black solid curve represents the limit cycle $x(t)$, the black dashed line represents the vertical line passing through the zero of
$A(t)$, and the blue dashed curve represents $\varphi(t)$.
The two subfigures illustrate the subcases where the limit cycle intersects
$\varphi(t)$ exactly once at $(t_{1},x_{1})$ when $A(t)$ has either a single simple zero or two simple zeros.
    }}\label{figure1}
\end{figure}

\noindent\textbf{ \textbf{Subcase 1.1}}: $Z(A)$ is a one-element set $\{a\}$.
We mainly focus on the case $a>t_{1}$, as illustrated in Fig. \ref{figure1}(A).
For the symmetric case $a<t_{1}$, the same argument applies and the same conclusion follows.
It is  {clear} that $x(t)-\varphi(t)$ is continuous and does not vanish on $(0,T)\backslash\{t_{1},a\}$.
Since $x(t_{1})-\varphi(t_{1})=0$ and $x'(t_{1})-\varphi'(t_{1})=-\varphi'(t_{1})<0$,
it follows that
\begin{align*}
\begin{split}
x(t)-\varphi(t)
\left\{
      \begin{aligned}
            &>0 ,& t\in(0,t_{1}), \\
            &<0 ,& t\in(t_{1},a). \\
        \end{aligned}
      \right.
\end{split}
\end{align*}
Moreover, since $\varphi(t)$ is increasing, $\lim_{t\rightarrow a^{+}}\varphi(t)=-\infty$. We get that
$$x(t)-\varphi(t)>0,\indent t\in(a,T).$$
Thus, by the definition of $\tau_{1}$ and $\tau_{2}$, these yield
\begin{equation}\label{eq121}
   x\left(\tau_{1}\left(s\right)\right)-\varphi\left(\tau_{1}\left(s\right)\right)>0,\indent s\in (x_{0},x_{1}),
\end{equation}
and
\begin{equation}\label{eq131}
    \begin{split}
        x\left(\tau_{2}\left(s\right)\right)-\varphi\left(\tau_{2}\left(s\right)\right)\left\{
      \begin{aligned}
            &>0 ,& s\in(x_{0},x(a)), \\
            &<0 ,& s\in(x(a),x_{1}). \\
        \end{aligned}
      \right.
    \end{split}
\end{equation}
We obtain by \eqref{eq011}, \eqref{eq121} and \eqref{eq131} that $W'(s)>0$ for $s\in(x(a), x_{1})$.
On the other hand, note that the expression in \eqref{case1W'} is well-defined on $(x_{0},x(a))$.
Also, by the monotonicity of $\varphi$, $\tau_{1}$, and $\tau_{2}$, one can easily verify that $\varphi\circ\tau_{1}-\varphi\circ\tau_{2}$ is strictly increasing on $(x_{0},x(a))$.
Therefore, $W'(s)$ has at most one sign change in $(x_{0},x_1)$.
Combining $W(x_{0})=W(x_{1})=0$, we get $W(s)<0$ for $s\in(x_{0},x_{1})$.

\noindent\textbf{ \textbf{Subcase 1.2}}: $Z(A)$ is a two-element set $\{a,b\}$ with $a<b$. See Fig. \ref{figure1}(B) for an illustration.
First, since $\varphi(t)$ is increasing,  $\lim_{t\rightarrow a^{+}}\varphi(t)=-\infty$ and $\lim_{t\rightarrow b^{-}}\varphi(t)=+\infty$. Therefore the unique zero $t_1$ of $x(t)-\varphi(t)$ must lie in the interval $(a,b)$.
Similar to the argument in Subcase 1.1, $x(t)-\varphi(t)$ is continuous and does not vanish on $(0,T)/\{a,t_{1},b\}$.
Recall that $x'(t_{1})-\varphi'(t_{1})=-\varphi'(t_{1})<0$.
We obtain that
\begin{align*}
\begin{split}
x(t)-\varphi(t)
\left\{
      \begin{aligned}
            &>0 ,& t\in(a,t_{1}), \\
            &<0 ,& t\in(t_{1},b). \\
        \end{aligned}
      \right.
\end{split}
\end{align*}
Moreover, the monotonicity of $\varphi$ also implies that $\lim_{t\rightarrow a^{-}}\varphi(t)=+\infty$ and $\lim_{t\rightarrow b^{+}}\varphi(t)=-\infty$. This yields
\begin{align*}
\begin{split}
x(t)-\varphi(t)
\left\{
      \begin{aligned}
            &<0 ,& t\in(0,a), \\
            &>0 ,& t\in(b,T). \\
        \end{aligned}
      \right.
\end{split}
\end{align*}
Hence, if we set $x_{*}=\min\{x(a),x(b)\}$ and $x^{*}=\max\{x(a),x(b)\}$, then from the definition of $\tau_{1}$ and $\tau_{2}$, the following inequalities hold.
\begin{equation}\label{eq231}
    \begin{split}
        x(\tau_{1}(s))-\varphi(\tau_{1}(s))\left\{
      \begin{aligned}
            &<0 ,& s\in(x_{0},x_{*}), \\
            &>0 ,& s\in(x^{*},x_{1}), \\
        \end{aligned}
      \right.
    \end{split}
\end{equation}
and
\begin{equation}\label{eq232}
    \begin{split}
        x(\tau_{2}(s))-\varphi(\tau_{2}(s))\left\{
      \begin{aligned}
            &>0 ,& s\in(x_{0},x_{*}), \\
            &<0 ,& s\in(x^{*},x_{1}). \\
        \end{aligned}
      \right.
    \end{split}
\end{equation}


Now according to \eqref{eq231} and \eqref{eq232}, we have that $W'(s)<0$ for $s\in(x_{0}, x_{*})$ and $W'(s)>0$ for $s\in(x^{*},x_{1})$.
We also emphasize that the expression in \eqref{case1W'} is well-defined on $(x_*,x^*)$. Moreover, the monotonicity of $\varphi$, $\tau_{1}$, and $\tau_{2}$ ensures that $\varphi\circ\tau_{1}-\varphi\circ\tau_{2}$ is increasing on $(x_{*},x^*)$.
As a consequence, $W'(s)$ can change sign at most once in $(x_{0},x_{1})$.
Taking into account $W(x_{0})=W(x_{1})=0$, we conclude that $W(s)<0$ for $s\in(x_{0},x_{1})$.
\vskip0.2cm

Based on the analyses in Subcase 1.1 and Subcase 1.2, we know by \eqref{express} that $P''(x_{0})>0$. Moreover, in the above omitted case $x_0>x_1$, a similar argument actually yields $W(s)>0$ for $s\in(x_{1},x_{0})$. Then, by \eqref{express} and noting the orientation of the integral (from $x_0$ to $x_1$), $P''(x_0)>0$ remains valid. The conclusion is verified for Case 1.
\vskip0.2cm

\noindent\textbf{ \textbf{Case 2}}: $x(t)-\varphi(t)$ has two different zeros in $(0,T)\backslash Z(A)$ (i.e., $x(t)$ has two extremum points), denoted by $t=t_{1}$ and $t=t_{2}$ with $t_{1}<t_{2}$. In this case, $x(t)$ is strictly monotonic on $[0,t_1]$, $[t_1,t_2]$ and $[t_2,T]$, respectively.
Again using Proposition \ref{YHLproposition}, if $x=x(t)$ is non-hyperbolic, then $h(T)=0$. Furthermore, the multiplicity and stability of $x=x(t)$ are given by
\begin{equation}\label{eq22-3}
    \begin{split}
        P&''(x_{0})\\
        &=-\frac{2}{x_{0}^{2}}\int_{0}^{T}\exp h(t)dx(t)\\
        &=-\frac{2}{x_{0}^{2}}\left(\int_{0}^{t_{1}}\exp h(t)dx(t)+\int_{t_1}^{t_{2}}\exp h(t)dx(t)+\int_{t_{2}}^{T}\exp h(t)dx(t)\right)\\
        &=-\frac{2}{x_{0}^{2}}\left(\int_{x_{0}}^{x_{1}}\exp h(\tau_{1}(x))dx+\int_{x_{1}}^{x_{2}}\exp h(\tau_{2}(x))dx+\int_{x_{2}}^{x_{0}}\exp h(\tau_{3}(x))dx\right),
    \end{split}
\end{equation}
where $x_{1}=x(t_{1})$, $x_{2}=x(t_{2})$, and $\tau_{1}$, $\tau_{2}$, $\tau_{3}$ represent the inverse functions of $x|_{[0,t_{1}]}$, $x|_{[t_{1},t_{2}]}$ and $x|_{[t_{2},T]}$, respectively.

Analogously to Case 1, we mainly consider the case $x_{2}<x_{0}<x_{1}$, and give a brief explanation for the opposite case $x_{2}>x_{0}>x_{1}$ in the end.
First, to simplify the expression in \eqref{eq22-3}, we introduce a function
\begin{equation}
    \begin{split}
        \tau_{0}(x)=
        \left\{
      \begin{aligned}
            &\tau_{3}(x),&  x\in [x_{2},x_{0}),\\
            &\tau_{1}(x),&  x\in [x_{0},x_{1}].\\
        \end{aligned}
      \right.
    \end{split}
    \end{equation}
It is clear that $\tau_0(x)$ is continuous on $[x_{2},x_{0})\cup(x_{0},x_{1}]$ and satisfies
\begin{align}\label{eq22-1}
  \lim_{x\rightarrow x_0^+}\tau_0(x)=\tau_1(x_0)=0,\ \ \
    \lim_{x\rightarrow x_0^-}\tau_0(x)=\tau_3(x_0)=T.
\end{align}
Moreover, we have
\begin{align}\label{eq22-2}
    &\tau_{0}'|_{(x_{2},x_{0})\cup(x_{0},x_{1})}>0,\  \tau_{2}'|_{(x_{2},x_{1})}<0,\  \tau_{0}(x_{1})=\tau_{2}(x_{1})=t_{1},\  \tau_{0}(x_{2})=\tau_{2}(x_{2})=t_{2}.
\end{align}

By the definition of $\tau_0$, \eqref{eq22-3} is reduced to
\begin{equation}\label{express1}
    P''(x_{0})=-\frac{2}{x_{0}^{2}}\left(\int_{x_2}^{x_1}\left(\exp h(\tau_{0}(x))-\exp h(\tau_{2}(x))\right)dx\right).
\end{equation}
The next step is to determine the sign of $P''(x_{0})$. For this purpose, we consider an auxiliary function
\begin{equation*}\label{eq22-4}
    W(s):=h(\tau_{0}(s))-h(\tau_{2}(s))=\int_{\tau_{2}(s)}^{\tau_{0}(s)}A(t)x^2(t)dt.
\end{equation*}
It is easy to know from \eqref{eq22-1} that
\begin{equation*}
    \lim_{s\rightarrow x_{0}^{-}}W(s)-\lim_{s\rightarrow x_{0}^{+}}W(s)=\int_{\tau_{2}(x_{0})}^{T}A(t)x^2(t)dt-\int_{\tau_{2}(x_{0})}^{0}A(t)x^2(t)dt=h(T)=0.
\end{equation*}
Thus, $W(s)$ is continuous on $(x_{2},x_{1})$ and differentiable on $(x_{2},x_{1})/\{x_{0}\}$. In addition, note that $W(x_{2})=W(x_{1})=0$ by \eqref{eq22-2}.
The problem becomes analyzing the sign of $W'(s)$. Clearly, $W'(s)$ is continuous on $(x_{2},x_{1})/\{x_{0}\}$, and one has
\begin{equation*}
    \begin{split}
        W'(s)&=A(\tau_{0}(s))x^2(\tau_{0}(s))\tau_{0}'(s)-A(\tau_{2}(s))x^2(\tau_{2}(s))\tau_{2}'(s)\\
        &=\frac{A(\tau_{0}(s))}{A(\tau_{0}(s))x(\tau_{0}(s))+B((\tau_{0}(s)))}-\frac{A(\tau_{2}(s))}{A(\tau_{2}(s))x(\tau_{2}(s))+B((\tau_{2}(s)))}.
    \end{split}
\end{equation*}
Recall that $Z(A)$ is the set of the simple zeros of $A(t)$ in $(0,T)$.
For $s\in (x_{2},x_{1})\setminus\big(\{x_{0}\}\cup \{x(t)|t\in Z(A)\}\big)$, we obtain
\begin{equation}\label{eq11}
    W'(s)=\frac{1}{x(\tau_{0}(s))-\varphi(\tau_{0}(s))}-\frac{1}{x(\tau_{2}(s))-\varphi(\tau_{2}(s))}.
\end{equation}
Since $x\circ\tau_{0}=x\circ\tau_{0}=\text{id}$, \eqref{eq11} can be written as
\begin{equation}\label{case2W'}
    W'(s)=\frac{\varphi(\tau_{0}(s))-\varphi(\tau_{2}(s))}{\left(x(\tau_{0}(s))-\varphi(\tau_{0}(s))\right)\left(x(\tau_{2}(s))-\varphi(\tau_{2}(s))\right)}.
\end{equation}
Then, the following discussion is subdivided into two subcases:

\begin{figure}[p]
    \centering
    \begin{subfigure}[b]{0.8\textwidth}
        \includegraphics[scale=0.7]{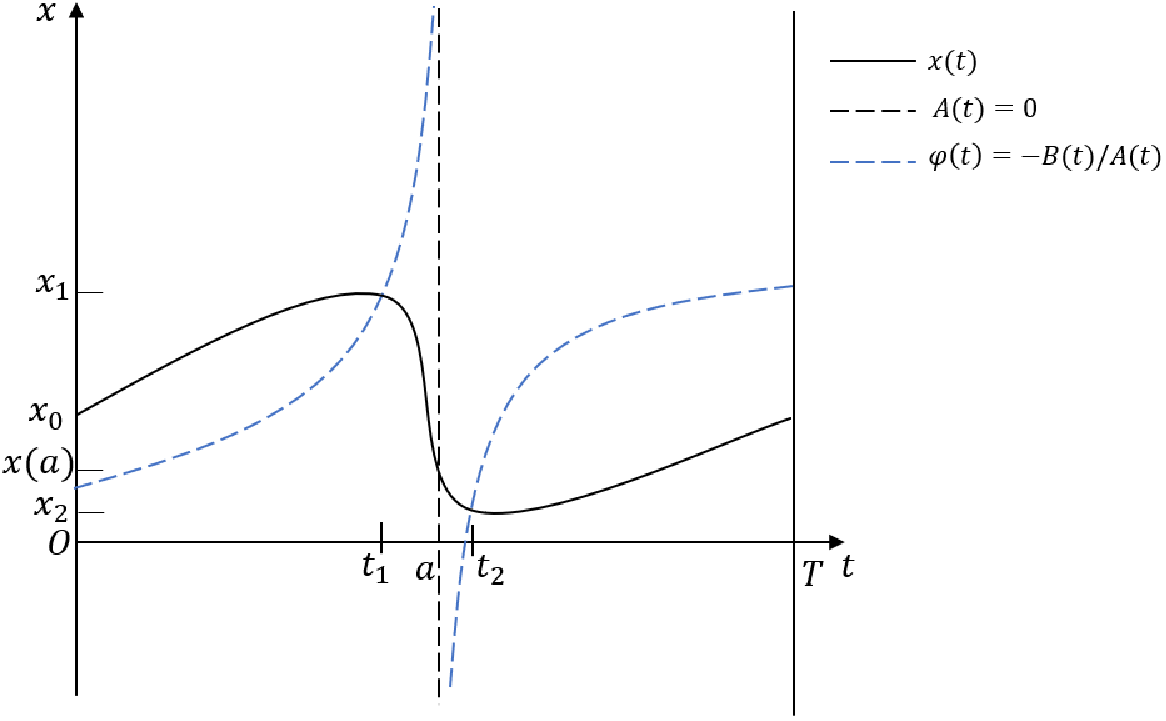}
        \caption{ {Subcase where $Z(A)$ is a one-element set $\{a\}$.}}
    \end{subfigure}
    \begin{subfigure}[b]{0.8\textwidth}
            \includegraphics[scale=0.7]{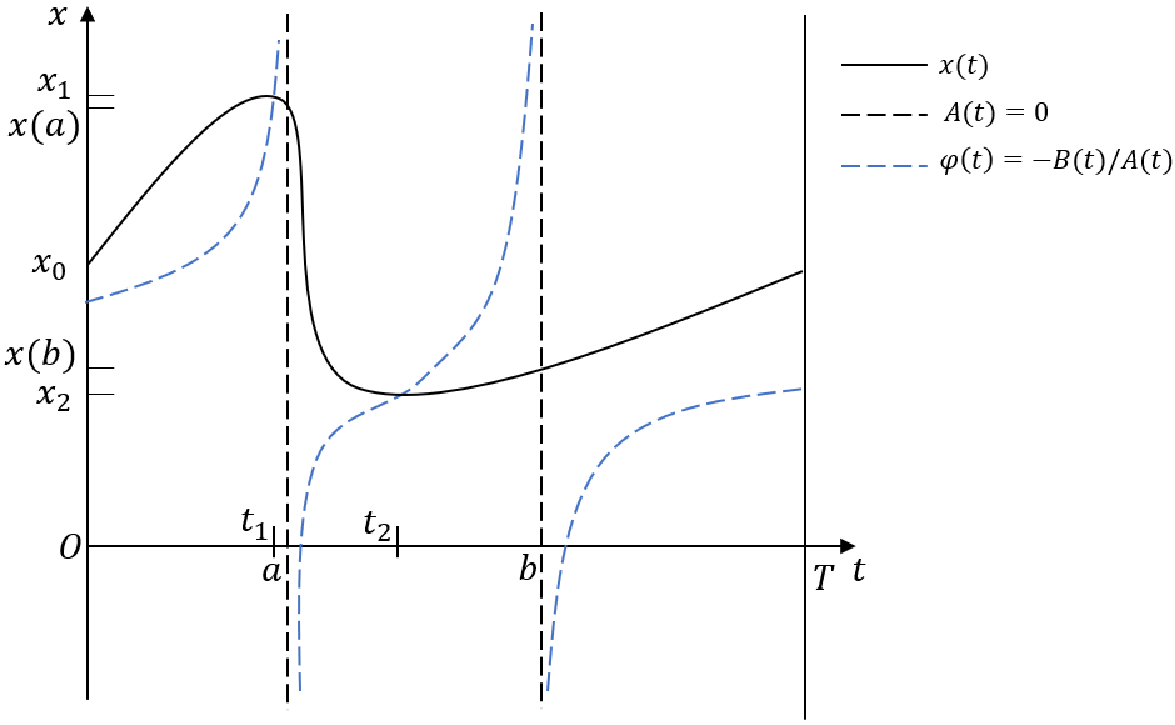}
            \caption{ {Subcase where $Z(A)$ is a two-element set $\{a,b\}$ with $a<b$.}}
    \end{subfigure}
    \caption{Case where $x(t)-\varphi(t)$ has exactly two zeros $t=t_{1}$ and $t=t_{2}$ in $[0,T)$.
    {The black solid curve represents the limit cycle $x(t)$, the black dashed line represents the vertical line passing through the zero of
$A(t)$, and the blue dashed curve represents $\varphi(t)$.
The two subfigures illustrate the subcases where the limit cycle intersects
$\varphi(t)$ twice at $(t_{1},x_{1})$ and $(t_{2},x_{2})$ with $t_{1}<t_{2}$ when $A(t)$ has either a single simple zero or two simple zeros.
    }
    }\label{figure2}
\end{figure}

\noindent\textbf{ \textbf{Subcase 2.1}}: $Z(A)$ is a one-element set $\{a\}$.
See Fig. \ref{figure2}(A) for an illustration.
According to statement (i) of Lemma \ref{lemma2}, the two zeros of $x(t)-\varphi(t)$ must be located in the intervals $(0,a)$ and $(a,T)$, respectively.
Thus, $0<t_{1}<a<t_{2}<T$. We know that $x(t)-\varphi(t)$ is continuous and does not vanish on $(0,T)\setminus\{a,t_{1},t_{2}\}$.
Since $x(t_{i})-\varphi(t_{i})=0$ and $x'(t_{i})-\varphi'(t_{i})=-\varphi'(t_{1})<0$ for $i=1,2$,
it follows that
\begin{align*}
\begin{split}
x(t)-\varphi(t)
\left\{
      \begin{aligned}
            &>0 ,& t\in(0,t_{1})\cup(a,t_2), \\
            &<0 ,& t\in(t_{1},a)\cup(t_2,T). \\
        \end{aligned}
      \right.
\end{split}
\end{align*}
This together with the monotonicity of $\tau_{0}$ and $\tau_{2}$ yields
\begin{equation}\label{eq421}
    \begin{split}
        x(\tau_{0}(s))-\varphi(\tau_{0}(s))\left\{
      \begin{aligned}
            &<0 ,& s\in(x_{2},x_{0}), \\
            &>0 ,& s\in(x_{0},x_{1}). \\
        \end{aligned}
      \right.
    \end{split}
\end{equation}
and
\begin{equation}\label{eq431}
    \begin{split}
        x(\tau_{2}(s))-\varphi(\tau_{2}(s))\left\{
      \begin{aligned}
            &>0 ,& s\in(x_{2},x(a)), \\
            &<0 ,& s\in(x(a),x_{1}). \\
        \end{aligned}
      \right.
    \end{split}
\end{equation}


Now define $x_{*}=\min\{x(a),x_{0}\}$ and $x^{*}=\max\{x(a),x_{0}\}$. It follows from \eqref{eq11}, \eqref{eq421} and \eqref{eq431} that $W'(s)<0$ for $s\in(x_2,x_*)$ and $W'(s)>0$ for $s\in(x^*,x_1)$.
On the other hand, note that the expression in \eqref{case2W'} applies on $(x_{*},x^{*})$. Also, by the monotonicity of $\varphi$, $\tau_{0}$ and $\tau_{2}$,
one can verify that $\varphi\circ\tau_{0}-\varphi\circ\tau_{2}$ is increasing on $(x_{*},x^{*})$.
Therefore, $W'(s)$ has at most one sign change in $(x_{*},x^{*})$, and by its continuity on $(x_{2},x_{1})\setminus\{x_{0}\}$ the same holds in $(x_{2},x_{1})$.
Combining $W(x_{0})=W(x_{1})=0$, we obtain that $W(s)<0$ in $(x_{2},x_{1})$.

\noindent\textbf{ \textbf{Subcase 2.2}}: $Z(A)$ is a two-element set $\{a,b\}$ with $a<b$.
From statement (i) of Lemma again, the zeros $t_1$ and $t_2$ of $x(t)-\varphi(t)$ lie in two distinct intervals among $(0,a)$, $(a,b)$, and $(b,T)$.
Moreover, due to $\lim_{t\rightarrow a^{+}}\varphi(t)=-\infty$ and $\lim_{t\rightarrow b^{-}}\varphi(t)=+\infty$,
one of $t_1$ and $t_2$ must belong to the interval $(a,b)$.
For these reasons, we assume without loss of generality that $0<t_{1}<a<t_{2}<b<T$, as depicted in Fig. \ref{figure2}(B).
The same conclusion for the symmetric case $0<a<t_{1}<b<t_{2}<T$ follows in exactly the same way.

Clearly, $x(t)-\varphi(t)$ is continuous and does not vanish on $ (0,T)\setminus\{t_{1},a,t_{2},b\}$.
Since $x(t_{i})-\varphi(t_{i})=0$ and $x'(t_{i})-\varphi'(t_{i})=-\varphi'(t_{i})<0$ for $i=1,2$, we have
\begin{align*}
\begin{split}
x(t)-\varphi(t)
\left\{
      \begin{aligned}
            &>0 ,& t\in(0,t_{1})\cup(a,t_2), \\
            &<0 ,& t\in(t_{1},a)\cup(t_2,b). \\
        \end{aligned}
      \right.
\end{split}
\end{align*}
Moreover, observe that $\lim_{t\rightarrow b^{+}}\varphi(t)=-\infty$. Therefore $x(t)-\varphi(t)>0$ for $t\in(b,T)$.
Taking into account the monotonicity of $\tau_{0}$ and $\tau_{2}$, these yield
\begin{equation}\label{eq521}
    \begin{split}
        x(\tau_{0}(s))-\varphi(\tau_{0}(s))\left\{
      \begin{aligned}
            &<0 , &s\in(x_{2},x(b)), \\
            &>0 , &s\in (x(b),x_{0})\cup(x_{0},x_{1}), \\
        \end{aligned}
      \right.
    \end{split}
\end{equation}
and
\begin{equation}\label{eq531}
    \begin{split}
        x(\tau_{2}(s))-\varphi(\tau_{2}(s))\left\{
      \begin{aligned}
            &>0 ,& s\in(x_{2},x(a)), \\
            &<0 ,& s\in(x(a),x_{1}). \\
        \end{aligned}
      \right.
    \end{split}
\end{equation}


Now let $x_{*}=\min\{x(a),x(b)\}$ and $x^{*}=\max\{x(a),x(b)\}$.
If $x_{0}\not\in(x_{*},x^{*})$, then
from \eqref{eq11}, \eqref{eq521} and \eqref{eq531}, we get $W'(s)<0$ for $s\in(x_{2},x_{*})$ and $W'(s)>0$ for $s\in(x^{*}, x_{1})\setminus\{x_{0}\}$.
Moreover, note that $W'(s)$ is continuous on $(x_{*},x^{*})$ and the expression in \eqref{case2W'} applies. Also, the monotonicity of $\varphi$, $\tau_{0}$ and $\tau_{2}$ ensures that $\varphi\circ\tau_{0}-\varphi\circ\tau_{2}$ is increasing on $(x_{*},x^{*})$.
Therefore, $W'(s)$ can change sign at most once in $(x_{*},x^{*})$, and by its continuity on $(x_{2},x_{1})\setminus\{x_{0}\}$ the same is true in $(x_{2},x_{1})$.
Combining $W(x_{2})=W(x_{1})=0$, we obtain that $W(s)<0$ in $(x_{2},x_{1})$.

If $x_{0}\in(x_{*},x^{*})$, then \eqref{eq11}, \eqref{eq521} and \eqref{eq531} indicate that $W'(s)<0$ for $s\in(x_{2},x_{*})$ and $W'(s)>0$ for $s\in(x^{*}, x_{1})$.
It remains to consider the sign of $W'(s)$ on $(x_*,x^*)\setminus\{x_0\}$. To this end, we first know by \eqref{eq521} and \eqref{eq531} that the denominator of the expression in \eqref{case2W'} is always positive on $(x_{*}, x_{0})\cup(x_{0},x^{*})$. Next, we show that $\varphi\circ\tau_{0}-\varphi\circ\tau_{2}$ is increasing on the whole $(x_{*},x^{*})$. In fact, due to the monotonicity of $\varphi$, $\tau_{0}$ and $\tau_{2}$, the function $\varphi\circ\tau_{0}-\varphi\circ\tau_{2}$ is increasing on both $(x_{*}, x_{0})$ and $[x_{0},x^{*})$. Observe that
\begin{equation*}
    \lim_{s\rightarrow x_{0}^{-}}\varphi \left(\tau_{0}\left(s\right)\right)-\varphi\left(\tau_{2}\left(s\right)\right)=\lim_{s\rightarrow x_{0}^{-}}\varphi \left(\tau_{3}\left(s\right)\right)-\varphi\left(\tau_{2}\left(s\right)\right)=\varphi(T)-\varphi\left(\tau_{2}\left(x_{0}\right)\right),
\end{equation*}
and
\begin{equation*}
    \varphi \left(\tau_{0}\left(x_0\right)\right)-\varphi\left(\tau_{2}\left(x_0\right)\right)=\varphi(0)-\varphi\left(\tau_{2}\left(x_{0}\right)\right).
\end{equation*}
Thus, it is sufficient to prove $\varphi(0)\geq \varphi(T)$. Here we assume for a contradiction that $\varphi(0)< \varphi(T)$.
Then, again using
\begin{align*}
  \lim_{t\rightarrow a^{+}}\varphi(t)=\lim_{t\rightarrow b^{+}}\varphi(t)=-\infty,\indent
  \lim_{t\rightarrow a^{-}}\varphi(t)=\lim_{t\rightarrow b^{-}}\varphi(t)=+\infty,
\end{align*}
there exists $\lambda\in (\varphi(0),\varphi(T))$ such that $\varphi(t)-\lambda=-\frac{\lambda A(t)+B(t)}{A(t)}$ has three zeros located in $(0,a)$, $(a,b)$, and $(b,T)$, respectively.
This contradicts the assumption that $\{A,B\}$ is an ET-system with accuracy one on $[0,T)$.
Hence, $\varphi(0)\geq \varphi(T)$ holds and therefore $\varphi\circ\tau_{0}-\varphi\circ\tau_{2}$ is increasing on $(x_{*},x^{*})$. As a consequence, $W'(s)$ can change sign at most once in $(x_*,x^*)\setminus\{x_0\}$, and by its continuity on $(x_2,x_1)\setminus\{x_0\}$ the same is true in $(x_2,x_1)$. Taking into account $W(x_{0})=W(x_{1})=0$, we get that $W(s)<0$ in $(x_{2},x_{1})$ .
\vskip0.2cm

According to the analyses in Subcase 2.1 and Subcase 2.2, we finally get by \eqref{express} that $P''(x_{0})>0$.
 We emphasize that in the omitted case $x_{2}>x_{0}>x_{1}$, a similar argument yields $W(s)>0$ for $s\in(x_{1},x_{0})$. Then, by \eqref{express} and noting the orientation of the integral (from $x_0$ to $x_1$), $P''(x_0)>0$ remains valid.
In summary, the conclusion of the statement holds for Case 2.

The proof is finished.
\end{proof}

\subsection{The Lyapunov constants of solution $x=0$}
This subsection studies the local dynamics of the  {trivial} solution $x=0$ of equation \eqref{eq0}.
For this purpose, let $x=x(t,x_0)$ be the solution of equation \eqref{eq0} with initial condition $x(0,x_0)=x_0$. It is clear that $x=x(t,x_0)$ is well-defined on $[0,T]$ for all $x_0$ in a small neighbourhood of zero. Moreover, due to the smoothness of the equation in variable $x$, $x(t,x_0)$ can be expanded into a power series as:
\begin{equation*}
    x(t,x_{0})=x_{0}+\sum_{i=2}^{+\infty}a_{i}(t)x_{0}^{i},\indent t\in[0,1].
\end{equation*}
Then, the displacement function $d(x_0)$ of the equation near $x_0=0$ is given by
\begin{equation*}
    d(x_{0}):=x(T,x_{0})-x_{0}=\sum_{i=2}^{+\infty}V_{i}x_{0}^{i},\indent V_{i}=a_{i}(T).
\end{equation*}
For each $i\in\mathbb{N}_{\geq 2}$, $V_i$ is called the \textit{$i$-th Lyapunov constant} of the trivial solution $x=0$ (see for instance \cite{YHL,rotated2,rotated3,AGG}).
We stress that when $V_i$ is the first non-zero Lyapunov constant, the solution $x=0$ forms a limit cycle with multiplicity $i$, and its stability is determined by the sign of $V_i$.

For equation \eqref{eq0}, a straightforward calculation or the direct application of the result in \cite{AL1987} yields
\begin{equation}\label{eq3.2}
    \begin{split}
        V_{2}&=\int_{0}^{T}B(t)dt,\\
        V_{3}&=\int_{0}^{T}A(t)dt, \text{ if }V_{2}=0,\\
        V_{4}&=\int_{0}^{T}A(t)\int_{0}^{t}B(s)dsdt, \text{ if } V_{2}=V_{3}=0.
    \end{split}
\end{equation}
Then, we have the following result.

\begin{proposition}\label{Lyapunov constants}
    Assume that the set of coefficients $\{A,B\}$ of equation \eqref{eq0} forms an ET-system with accuracy one on the interval $[0,T)$.
    Then, one of the following statements holds.
    \begin{itemize}
        \item[(i)] The set $\{A,B\}$ satisfies $\bm{(D.1)}$ with $E=[0,T)$, and $x=0$ is a center (resp. the unique limit cycle) if the Lyapunov constants $V_{2}=V_{3}=0$ (resp. $V_2=0$ and $V_3\neq0$).
        \item[(ii)] The set $\{A,B\}$ satisfies $\bm{(D.2)}$ with $E=[0,T)$, and $x=0$ is a limit cycle with multiplicity at most four.
        Moreover, if the Lyapunov constants $V_{2}=V_{3}=0$ and $A'(t)B(t)-A(t)B'(t)>0$ (resp. $<0$), then $V_{4} > 0$ (resp. $< 0$).
    \end{itemize}
\end{proposition}
\begin{proof}
    According to the quasi-dichotomy in Proposition \ref{lemma1}, it is sufficient to prove the assertion for $x=0$ in each statement.

    (i) We begin by showing that $x=0$ is a center if $V_{2}=V_{3}=0$.
        In this case, we have $\int_{0}^{T}A(t)\,dt=\int_{0}^{T}B(t)\,dt=0$ from \eqref{eq3.2}.
        Then, for $(\lambda,\mu)\neq(0,0)$ given in $\bm{(D.1)}$, the function $\lambda A(t)+\mu B(t)$ simultaneously does not change sign on $[0,T)$ and satisfies
        \begin{align*}
            \int_{0}^{T}\big(\lambda A(t)+\mu B(t)\big)\,dt=0.
        \end{align*}
        This implies that $\lambda A(t)+\mu B(t)\equiv 0$. Therefore, equation \eqref{eq0} becomes either $\frac{dx}{dt}=B(t)\left(-\frac{\mu}{\lambda}x^{3}+x^{2}\right)$ when $\lambda\not=0$, or $\frac{dx}{dt}=A(t)\left(x^{3}-\frac{\lambda}{\mu}x^{2}\right)$ when $\mu\not=0$.
        In both cases all the solutions of the equation are closed.
        Consequently, $x=0$ is a center if $V_{2}=V_{3}=0$.

        To prove that $x=0$ is the unique limit cycle when $V_2=0$ and $V_3\neq0$, we divide the argument into the following {three cases}.

    \noindent{ Case 1}: $\lambda=0$. We have $B(t)\geq 0$ (or $\leq 0$).
    Since $V_{2}=0$, it follows from \eqref{eq3.2} that $B(t)\equiv 0$.
    Thus, equation \eqref{eq0} is reduced to $\dot{x}=A(t)x^{3}$, which clearly has no limit cycles in the region $x\not=0$.

    \noindent{ Case 2}: $\lambda\not=0$ {and $\mu\not=0$}.
    Equation \eqref{eq0} can be written as
    \begin{equation}\label{890}
        {\frac{dx}{dt}=-\frac{\mu}{\lambda}B(t)x^{2}\left(x-\frac{\lambda}{\mu}\right)+\frac{1}{\lambda}\left(\lambda A(t)+\mu B(t)\right)x^{3},}
    \end{equation}
    where $\frac{1}{\lambda}\left(\lambda A(t)+\mu B(t)\right)x^{3}\geq0$ (or $\leq0$) on both regions $x>0$ and $x<0$.
    Note that $\int_{0}^{T} B(t)\,dt=V_2=0$. Thus, $x=0$ is a center of the equation
\begin{equation}\label{098}
    \frac{dx}{dt}=-\frac{\mu}{\lambda}B(t)x^{2}(x-\lambda).
\end{equation}
Comparing equation \eqref{890} with \eqref{098}, we know that equation \eqref{eq0} has no limit cycles in the region $x\not=0$.

    \noindent{ {Case 3}}: {
        $\lambda\not=0$ and $\mu=0$}.
{We have $A(t)\geq 0$ (or $\leq 0$).
Since $V_{2}=0$, it follow that $x=0$ is a center of the equation
\begin{equation}\label{B(t)}
    \frac{dx}{dt}=B(t)x^{2}.
\end{equation}
Comparing equation \eqref{eq0} with \eqref{B(t)}, we know that equation \eqref{eq0} has no limit cycles in the region $x\not=0$.
Statement (i) is verified.}

    (ii) We only need to determine the Lyapunov constant $V_4$ of $x=0$ in the case $V_2=V_3=0$. Furthermore, taking the change of variable $x\rightarrow -x$ if necessary, we may suppose without loss of generality that $A'(t)B(t)-A(t)B'(t)>0$ on $[0,T)$. In the following we show that $V_4>0$.

    To this end, denote by $\mathcal{V}_{2}(t)=\int_{0}^{t}B(s)\,ds$ and $\mathcal{V}_{3}(t)=\int_{0}^{t}A(s)\,ds$. We introduce two auxiliary functions
    \begin{align}\label{def1}
      f(t)=\mu_{1}A(t)-\mu_{2}B(t),\indent G(t)=\lambda_{1}\mathcal{V}_{3}(t)-\lambda_{2}\mathcal{V}_{2}(t),
    \end{align}
    where $\mu_1, \mu_2, \lambda_1, \lambda_2\in\mathbb R$ are parameters to be chosen later. It is easy to get by \eqref{eq3.2} that
    \begin{equation}\label{ert'}
        \mathcal{V}_{2}(0)=\mathcal{V}_{2}(T)=V_{2}=0,\indent \mathcal{V}_{3}(0)=\mathcal{V}_{3}(T)=V_{3}=0,
    \end{equation}
    and
    \begin{equation*}\label{ert}
        V_{4}=\int_{0}^{T}\mathcal{V}_{2}(t)\,d\mathcal{V}_{3}(t)=-\int_{0}^{T}\mathcal{V}_{3}(t)\,d\mathcal{V}_{2}(t)=-\int_{0}^{T}B(t)\int_{0}^{t}A(s)\,dsdt.
    \end{equation*}
    Then, one has
    \begin{equation*}
    \begin{split}
        \int_{0}^{T}f(t)G(t)dt=&\int_{0}^{T}\left(\mu_{1}A(t)-\mu_{2}B(t)\right)\left(\lambda_{1}V_{3}(t)-\lambda_{2}V_{2}(t)\right)dt\\
        =&\int_{0}^{T}\left(\mu_{1}A(t)-\mu_{2}B(t)\right)\int_{0}^{t}\left(\lambda_{1}A(s)-\lambda_{2}B(s)\right)dsdt\\
        =&\frac{1}{2}\mu_{1}\lambda_{1} V_{3}^{2}+\frac{1}{2}\mu_{2}\lambda_{2} V_{2}^{2}-\mu_{1}\lambda_{2}V_{4}+\mu_{2}\lambda_{1}V_{4}\\
        =&(\mu_{2}\lambda_{1}-\mu_{1}\lambda_{2})V_{4}.
    \end{split}
\end{equation*}
Therefore, our argument is reduced to showing that $(\mu_{2}\lambda_{1}-\mu_{1}\lambda_{2})\int_{0}^{T}f(t)G(t)dt>0$ holds for an appropriate choice of $\mu_1, \mu_2, \lambda_1, \lambda_2$.

    Now let us begin by taking $\mu_{1}=B(0)$ and $\mu_{2}=A(0)$. We have
    \begin{align}\label{eq3.23}
      f(0)=0,\indent f'(0)=A'(0)B(0)-A(0)B'(0)> 0.
    \end{align}
    Then, $f$ is a non-trivial linear combination of $A$ and $B$. Moreover, since $\int_{0}^{T}f(t)\,dt=\mu_{1}V_{3}-\mu_{2}V_{2}=0$, there exists $t_{0}\in(0,T)$ such that $f(t_{0})=0$. According to the assumption that $\{A,B\}$ forms an ET-system with accuracy one on $[0,T)$, $f(t)$ has exactly two zeros $t=0$ and $t=t_{0}$ in the interval, which are both simple.
    Taking into account \eqref{eq3.23} again, we get
    \begin{equation}\label{eq3.21}
    \begin{split}
       f(t)\left\{
      \begin{aligned}
            &>0 , &t\in(0,t_{0}), \\
            &<0 , &t\in (t_{0},T). \\
        \end{aligned}
      \right.
    \end{split}
\end{equation}

    Next, set $\lambda_{1}=\mathcal{V}_{2}(t_{0})$ and $\lambda_{2}=\mathcal{V}_{3}(t_{0})$.
    It follows from \eqref{def1} and \eqref{ert'} that $G(0)=G(t_{0})=G(T)=0$. Since $G'(t)=\lambda_{1}A(t)-\lambda_{2}B(t)$, it has at most two zeros in $[0,T)$, analogous to $f(t)$. Thus, $t=0$, $t=t_{0}$ and $t=T$ are the only three zeros of $G$, and they are all simple. Observe that
    \begin{equation*}
    G'(0)=\lambda_1A(0)-\lambda_2B(0)
    =\int_{0}^{t_{0}}B(s)A(0)-A(s)B(0)\,ds=-\int_{0}^{t_{0}}f(s)\,ds<0.
    \end{equation*}
    This yields
    \begin{equation}\label{eq3.22}
    \begin{split}
       G(t)\left\{
      \begin{aligned}
            &<0 , &t\in(0,t_{0}), \\
            &>0 , &t\in (t_{0},T). \\
        \end{aligned}
      \right.
    \end{split}
    \end{equation}

    Finally, utilizing \eqref{eq3.21} and \eqref{eq3.22}, we get that $\int_{0}^{T}f(t)G(t)dt<0$. Moreover, the inequality \eqref{eq3.21} also implies that
    \begin{equation*}
        \mu_{2}\lambda_{1}-\mu_{1}\lambda_{2}=-\int_{0}^{t_{0}}\mu_{1}A(s)-\mu_{2}B(s)ds=-\int_{0}^{t_{0}}f(s)\,ds<0.
    \end{equation*}
    Consequently, $(\mu_{2}\lambda_{1}-\mu_{1}\lambda_{2})\int_{0}^{T}f(t)G(t)dt>0$ holds, which immediately leads to $V_4>0$. Statement (ii) follows.

    The proof is finished.
\end{proof}

\section{Maximum number of limit cycles of equation \eqref{eq0}}
The goal of this section is to prove Theorem \ref{theorem0}, utilizing Proposition \ref{multiplicity of non-zero limit cycle}, Proposition \ref{Lyapunov constants} and the theory of rotated equations. We again emphasize that the set of coefficients $\{A,B\}$ of equation \eqref{eq0} in the theorem forms an ET-system with accuracy one on $[0,T)$. Then according to statement (i) of Proposition \ref{multiplicity of non-zero limit cycle}, it is sufficient to consider the case where $\{A,B\}$ satisfies $\bm{(D.2)}$ but not $\bm{(D.1)}$, with $E=[0,T)$. Furthermore, by performing the change of variable $x\rightarrow -x$ if necessary, we may assume without loss of generality that $A'(t)B(t)-A(t)B'(t)>0$ on $[0,T)$. In summary, we only need to focus on equation \eqref{eq0} under the assumption of Theorem \ref{theorem0} with the following additional hypothesis:
\begin{itemize}
    \item[$\bm{(H)}$] $A'(t)B(t)-A(t)B'(t)>0$ on $[0,T)$, and $\bm{(D.1)}$ with $E=[0,T)$ is not satisfied.
\end{itemize}


In what follows we divide the argument into several steps.
It is clear that the coefficients $A,B$ of equation \eqref{eq0} in Theorem \ref{theorem0} can be written as
\begin{align}\label{A-B}
  A=\lambda_{0}f_{0}+\lambda_{1}f_{1}+\lambda_{2}f_{2},\indent
  B=\mu_{0}f_{0}+\mu_{1}f_{1}+\mu_{2}f_{2},
\end{align}
where $f_0,f_1,f_2$ are defined {as} in the theorem and $\lambda_{0},\lambda_{1},\lambda_{2},\mu_{0},\mu_{1},\mu_{2}\in\mathbb R$. Thus, equation \eqref{eq0} becomes
\begin{equation}\label{equation1}
    \begin{split}
        \frac{dx}{dt}
        =(\lambda_{0}f_{0}+\lambda_{1}f_{1}+\lambda_{2}f_{2})x^{3}+(\mu_{0}f_{0}+\mu_{1}f_{1}+\mu_{2}f_{2})x^{2}.
    \end{split}
\end{equation}
Let $\bm{\eta}=(\lambda_{0},\lambda_{1},\lambda_{2},\mu_{0},\mu_{1},\mu_{2})\in\mathbb{R}^{6}$ be the parameter vector associated with the equation.
Let $\mathcal{L}_{H}$ and $\mathcal{L}_{1}$ be the sets of $\bm\eta$ given by
\begin{align*}
\begin{split}
    &\mathcal{L}_{H}=\left\{\left.\bm{\eta}\in\mathbb{R}^{6}\right|\{A,B\} \text{ satisfies $\bm{(H)}$}\right\},\\
    &\mathcal{L}_{1}=\left\{\left.\bm{\eta}\in\mathbb{R}^{6}\right|\{ A,B\} \text{ satisfies } \bm{(D.1)} \text{ with } E=[0,T)\right\}.
\end{split}
\end{align*}
In the first step, we present the following characterization for $\mathcal L_H$, which will be useful later when applying the theory of rotated equations.

\begin{lemma}\label{LHopen}
    $\mathcal{L}_{H}$ is an open set in $\mathbb{R}^{6}$.
    Moreover, its boundary $\partial\mathcal{L}_{H}\subseteq\mathcal{L}_{1}$.
\end{lemma}
\begin{proof}
    Let $A$ and $B$ be given in \eqref{A-B}. To verify the assertion, we additionally introduce two sets
\begin{align*}
    &\mathcal{L}_{2}^{+}=\left\{\bm{\eta}\in\mathbb{R}^{6}\left|\ A'(t)B(t)-A(t)B'(t)>0 \text{ on } [0,T)\right.\right\},\\   &\mathcal{L}_{2}^{-}=\left\{\bm{\eta}\in\mathbb{R}^{6}\left|\ A'(t)B(t)-A(t)B'(t)<0 \text{ on } [0,T)\right.\right\}.
\end{align*}
Clearly, $\mathcal{L}_{2}^{+}\cup\mathcal{L}_{2}^{-}=\left\{\left.\bm{\eta}\in\mathbb{R}^{6}\right|\{A,B\} \text{ satisfies } \bm{(D.2)}\text{ with } E=[0,T)\right\}$. It follows from Proposition \ref{lemma1} that $\mathbb{R}^{6}=\mathcal{L}_{1}\cup \mathcal{L}_{2}^{+}\cup\mathcal{L}_{2}^{-}$. Also, by definition we have $\mathcal{L}_{H}=\mathcal{L}_{1}^{c}\cap\mathcal{L}_{2}^{+}$, where the superscript
``$c$'' denotes the complement of a set.

We first claim that $\mathcal{L}_{1}$ is closed in $\mathbb{R}^{6}$.
Denote by $A_{\bm\eta}(t)$ and $B_{\bm\eta}(t)$ the functions $A$ and $B$ associated with the parameter vector $\bm\eta$, respectively.
It is easy to see that
$$\mathcal{L}_{1}=\left\{\bm{\eta}\in\mathbb{R}^{6}\left|\ \exists (\lambda,\mu)\in \mathbb{S}^{1} \text{ s.t. } \lambda A_{\bm\eta}(t)+\mu B_{\bm\eta}(t)\geq 0 \text{ on } [0,T)\right.\right\}.$$
Suppose that $\bm{\eta}_{0}$ is a limit point of $\mathcal{L}_{1}$. Then there exists a sequence $\{\bm{\eta}_{n}\}_{n=1}^{+\infty}\subseteq \mathcal{L}_{1}$ such that $\bm{\eta}_{n}\rightarrow \bm{\eta}_{0}$ as $n\rightarrow +\infty$.
    For any $t\in[0,T)$, we have $A_{\bm{\eta}_{n}}(t)\rightarrow A_{\bm{\eta}_{0}}(t)$ and $B_{\bm{\eta}_{n}}(t)\rightarrow B_{\bm{\eta}_{0}}(t)$ as $n\rightarrow +\infty$.
    In addition, note that each $\bm{\eta}_{n}$ admits a corresponding pair of parameters $(\lambda_{n}, \mu_{n})\in \mathbb{S}^{1}$, satisfying $\lambda_{n}A_{\bm{\eta}_{n}}(t)+\mu_{n}B_{\bm{\eta}_{n}}(t)\geq 0$ on $[0,T)$.
    Due to the compactness of $\mathbb{S}^{1}$, the sequence $\{(\lambda_{n}, \mu_{n})\}_{n=1}^{+\infty}$ has a convergent subsequence $\{(\lambda_{n_{i}}, \mu_{n_{i}})\}_{i=1}^{+\infty}$ that converges to a point $(\lambda_{0},\mu_{0})\in\mathbb{S}^{1}$.
    Thus, we have $\lambda_{n_{i}}A_{\bm{\eta}_{n_{i}}}(t)+\mu_{n_{i}}B_{\bm{\eta}_{n_{i}}}(t)\rightarrow \lambda_{0}A_{\bm{\eta}_{0}}(t)+\mu_{0}B_{\bm{\eta}_{0}}(t)$ as $i\rightarrow +\infty$.
    This implies that
$\lambda_{0}A_{\bm{\eta}_{0}}(t)+\mu_{0}B_{\bm{\eta}_{0}}(t)\geq 0$ on $[0,T)$, and therefore $\bm{\eta}_{0}\in \mathcal{L}_{1}$.
    Hence, $\mathcal{L}_{1}$ is  {a} closed set in $\mathbb{R}^{6}$.

    Next, we show that $\partial\mathcal{L}_{H}\subseteq \mathcal{L}_{1}$.
    Observe that $\partial\mathcal{L}_{H}=\big(\partial\mathcal{L}_{H}\cap \overline{\mathcal{L}_{2}^{-}}\big)\cup \big(\partial\mathcal{L}_{H}\cap (\overline{\mathcal{L}_{2}^{-}})^{c}\big)$,
    where the notation `` $\overline{\cdot}$ '' represents the closure of a set.
    Thus, on the one hand,
    \begin{equation*}
        \partial\mathcal{L}_{H}\cap \overline{\mathcal{L}_{2}^{-}}\subseteq \overline{\mathcal{L}_{H}}\cap \overline{\mathcal{L}_{2}^{-}}\subseteq \overline{\mathcal{L}_{2}^{+}}\cap \overline{\mathcal{L}_{2}^{-}}
        =\left\{\bm{\eta}\left|A_{\bm{\eta}}'(t)B_{\bm{\eta}}(t)-A_{\bm{\eta}}(t)B_{\bm{\eta}}'(t)\equiv 0\right.\right\}\subseteq \mathcal{L}_{1}.
    \end{equation*}
On the other hand,
suppose that $\bm{\eta}\in \partial\mathcal{L}_{H}\cap (\overline{\mathcal{L}_{2}^{-}})^{c}$. Then there exists $\varepsilon_{0}>0$ such that for any $\varepsilon\in(0,\varepsilon_{0})$, the $\varepsilon$-neighborhood of $\bm\eta$ satisfies
$U_{\varepsilon}(\bm{\eta})\subset (\overline{\mathcal{L}_{2}^{-}})^c$.
Recall that $\mathbb{R}^{6}=\mathcal{L}_{1}\cup \mathcal{L}_{2}^{+}\cup\mathcal{L}_{2}^{-}$ and $\mathcal{L}_{H}=\mathcal{L}_{1}^{c}\cap\mathcal{L}_{2}^{+}$.
One has
$$U_{\varepsilon}(\bm{\eta})
\subseteq\mathcal{L}_{1}\cup\mathcal{L}_{2}^{+}=\mathcal{L}_{1}\cup\mathcal{L}_{H}.$$
Since $\bm{\eta}$ is a boundary point of $\mathcal{L}_{H}$ but not an interior point,
$U_{\varepsilon}(\bm{\eta})\cap\mathcal{L}_{1}\not= \emptyset$ always holds.
By the arbitrariness of $\varepsilon$ and the closedness of $\mathcal{L}_{1}$, we get $\bm{\eta}\in \mathcal{L}_{1}$, which yields that $\partial\mathcal{L}_{H}\cap (\overline{\mathcal{L}_{2}^{-}})^{c}\subseteq\mathcal L_1$.
Consequently, $\partial\mathcal{L}_{H}\subseteq \mathcal{L}_{1}$.

The above argument also indicates that
$$\partial\mathcal{L}_{H}\cap\mathcal{L}_{H}\subseteq\mathcal{L}_{1}\cap\big(\mathcal{L}_{1}^{c}\cap\mathcal{L}_{2}^{+}\big)=\emptyset.$$
Thus, $\mathcal{L}_{H}$ is an open set and the proof is finished.
\end{proof}





In the next step, we establish a preliminary result on the number of limit cycles of equation \eqref{equation1}. It is presented in the following subsection.

\subsection{Preliminary estimates for the number of limit cycles}
According to \eqref{eq3.2} and \eqref{A-B}, the second, third and fourth Lyapunov constants of $x=0$ for equation \eqref{equation1} are given by
\begin{align}\label{Lyapunov constants-2}
\begin{split}
    &V_{2}=\sum_{i=0}^{2}\mu_i I_i,\\
    &V_{3}=\sum_{i=0}^{2}\lambda_i I_i, \text{  if  } V_{2}=0, \\
    &V_{4}
    =\sum_{i,j=0}^{2}\lambda_i\mu_j I_{i,j},  \text{  if  } V_{2}=V_{3}=0,
\end{split}
\end{align}
where
\begin{align}\label{Lyapunov constants-2'}
\begin{split}
    I_i=\int_{0}^{T}f_{i}(t)\,dt,\ \
    I_{i,j}=\int_{0}^{T}f_{i}(t)\int_{0}^{t}f_{j}(s)\,ds\,dt,\indent i,j=0,1,2.
\end{split}
\end{align}
Let us begin by considering equation \eqref{equation1} under hypothesis $\bm{(H)}$.
From statement (ii) of Proposition \ref{Lyapunov constants}, the maximum multiplicity of $x=0$ is four.
Thus, there are exactly three cases for the equation: 
\begin{itemize}
  \item $V_{2}\not= 0$.
  \item $V_{2}= 0$ and $V_{3}\not= 0$.
  \item $V_{2}=V_{3}= 0$ and $V_{4}\not=0$.
\end{itemize}
We have the following result.

\begin{lemma}\label{V2V3}
    Suppose that the set of coefficients of equation \eqref{equation1} (i.e., equation \eqref{eq0} considered in Theorem \ref{theorem0}) satisfies hypothesis $\bm{(H)}$.
    Then one of the following statements holds.
    \begin{itemize}
        \item[(i)] $V_{2}\not= 0$.
        \begin{itemize}
            \item[(i.1)] If $V_{2}>0$, then equation \eqref{equation1} has at most two limit cycles (counted with multiplicities) in $x>0$ and $x<0$, respectively.
            \item[(i.2)] If $V_{2}<0$, then equation \eqref{equation1} has at most one limit cycle (counted with multiplicities) in $x>0$ and $x<0$, respectively.
        \end{itemize}
        \item[(ii)] $V_{2}=0$ and $V_{3}\not=0$.
        Then, equation \eqref{equation1} has at most one limit cycle (counted with multiplicity) in $x>0$ and $x<0$, respectively.
        \item[(iii)] $V_{2}=V_{3}= 0$ and $V_{4}\not=0$. Then, equation \eqref{equation1} has no limit cycles in the region $x\not= 0$.
    \end{itemize}
\end{lemma}
\begin{proof}
From the above argument, it is sufficient to verify the estimate for the number of limit cycles in the statements.
Furthermore, we may restrict our attention to the positive limit cycles because the conclusions for the negative limit cycles follow exactly in the same way.

For clarity, we introduce two variable parameters $\lambda ,\mu$ instead of $\lambda_0,\mu_0$, and then rewrite equation \eqref{equation1} as
\begin{align}\label{equation2}
 \frac{dx}{dt}=S(t,x;\lambda,\mu)
 =(\lambda f_{0}+\lambda_{1}f_{1}+\lambda_{2}f_{2})x^{3}+(\mu f_{0}+\mu_{1}f_{1}+\mu_{2}f_{2})x^{2}.
\end{align}
The associated parameter vector is denoted by $\eta_{\lambda,\mu}=(\lambda,\lambda_{1},\lambda_{2},\mu,\mu_{1},\mu_{2})$.  We also correspondingly denote the second, third and fourth Lyapunov constants in \eqref{Lyapunov constants-2} by $V_2(\mu)$, $V_3(\lambda)$ and $V_4(\lambda,\mu)$.
From assumption of Theorem \ref{theorem0}, we have $\frac{\partial S}{\partial \lambda}=f_0x^3>0$ and $\frac{\partial S}{\partial\mu}=f_0x^2>0$ for $x>0$. Thus, it follows from Definition \ref{definition of rotated} that \eqref{equation2} with fixed $\mu$ (resp. fixed $\lambda$) forms a family of rotated equations in $\lambda$ (resp. $\mu$) on $[0,T]\times\mathbb R^+$.

    (i.1) Let us take $\mu=\mu_0$ and consider the $\lambda$-parametric family \eqref{equation2}$|_{\mu=\mu_0}$. By assumption and Lemma \ref{LHopen}, there exists a maximal open interval $(\lambda_*,\lambda^*)$ containing $\lambda_0$, such that $\eta_{\lambda,\mu_0}\in\mathcal L_H$
    for $\lambda\in(\lambda_*,\lambda^*)$, where $\lambda_*,\lambda^*\in\mathbb R\cup\{-\infty,+\infty\}$. Observe that $\lambda f_{0}+\lambda_{1}f_{1}+\lambda_{2}f_{2}\leq 0$ on $[0,T]$ for
    \begin{equation*}
        \lambda\leq\lambda_{\text{min}}:=\min_{t\in[0,T]}\left\{-\frac{\lambda_{1}f_{1}+\lambda_{2}f_{2}}{f_{0}}\right\}.
    \end{equation*}
    Hence, $\eta_{\lambda,\mu_0}\in\mathcal L_1\subset \mathcal L^c_H$ for $\lambda\leq\lambda_{\text{min}}$.
    This implies $\lambda_{\text{min}}\leq\lambda_*$ and therefore $\lambda_*$ is finite.


    Assume for a contradiction that equation \eqref{equation2}$|_{\lambda=\lambda_0,\mu=\mu_0}$ (i.e. equation \eqref{equation1}) has $m\geq3$  positive limit cycles, counted with multiplicities.
    Without loss of generality, we may assume that all these limit cycles are hyperbolic.
    Indeed, according to hypothesis $\bm{(H)}$ and statement (ii) of Proposition \ref{multiplicity of non-zero limit cycle}, all the non-hyperbolic positive limit cycles of the equation must be lower-stable and upper-unstable.
    By Proposition \ref{property of rotated}, they split simultaneously in a small decrease of $\lambda$, each into two hyperbolic positive limit cycles.
    Then, equation \eqref{equation2}$|_{\lambda=\lambda_{0}-\varepsilon,\mu=\mu_0}$ has only hyperbolic positive limit cycles as $0<\varepsilon\ll 1$, and their number is at least $m$. By additionally choosing $\varepsilon>0$ sufficiently small such that $\lambda_{0}-\varepsilon\in(\lambda_*,\lambda^*)$, this becomes the case under consideration.

    Now we can denote by $x=x_{1}(t;\lambda)$, $x=x_{2}(t;\lambda)$, $\ldots$, $x=x_{m}(t;\lambda)$ the $m$ positive limit cycles of equation \eqref{equation2}$|_{\mu=\mu_0}$, with $\lambda$ being in some neighborhood of $\lambda_0$ and $0<x_{1}(t;\lambda_0)<\cdots<x_{m}(t;\lambda_0)$.
    Since $V_{2}(\mu_{0})>0$ from assumption, we know that $x=0$ is upper-unstable for equation \eqref{equation2}$|_{\lambda=\lambda_0,\mu=\mu_0}$.
    Thus, both $x=x_{1}(t;\lambda_{0})$ and $x=x_{3}(t;\lambda_{0})$ are stable, whereas $x=x_{2}(t;\lambda_{0})$ is unstable.
    According to Proposition \ref{property of rotated}, as $\lambda$ decreases from $\lambda_0$, $x=x_{1}(t;\lambda)$ and $x=x_{3}(t;\lambda)$ decrease whereas $x=x_{2}(t;\lambda)$ increases, unless a change of stability appears, as illustrated in Fig. \ref{figure3}.
    \begin{figure}[t]
    \centering
        \includegraphics[scale=0.45]{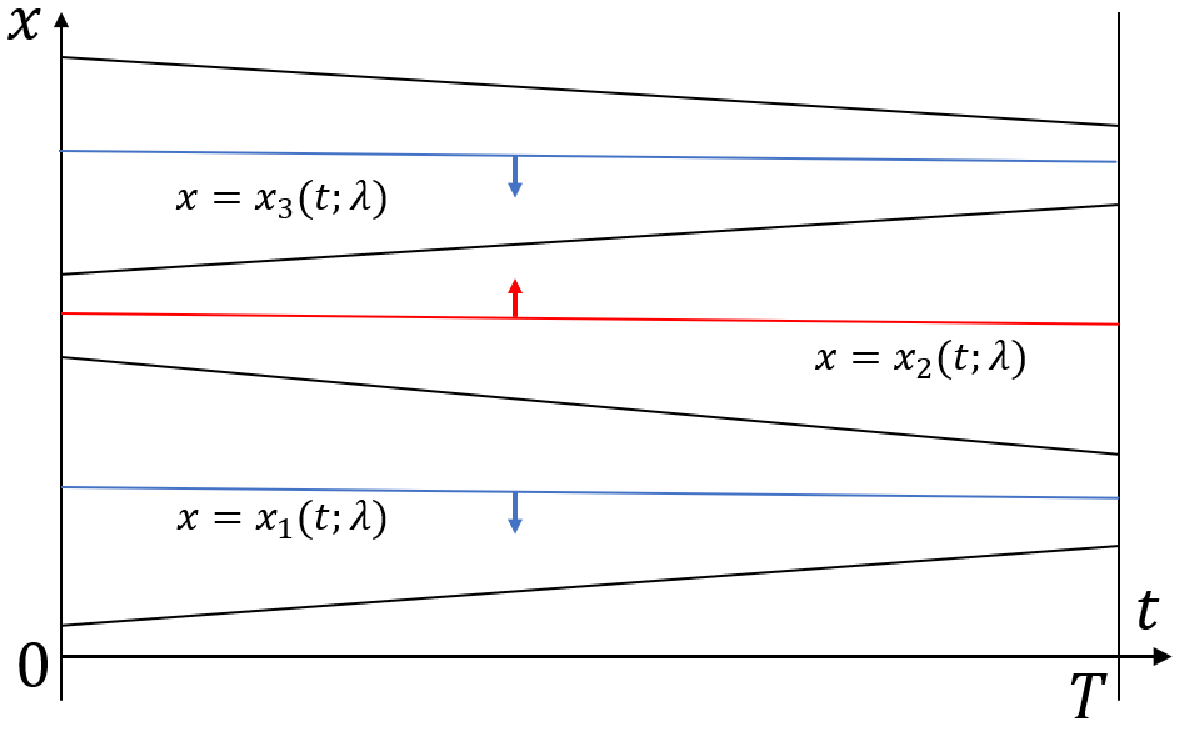}
    \caption{Direction of the movement of the stable and unstable limit cycles as parameter
$\lambda$ decreases. The blue, red and black lines represent stable limit cycles, unstable limit cycles and nearby non-periodic solutions, respectively.}\label{figure3}
\end{figure}

        Recall that when $\eta_{\lambda,\mu_0}\in\mathcal L_H$, statement (ii) of Proposition \ref{multiplicity of non-zero limit cycle} ensures that all the non-hyperbolic positive limit cycles of \eqref{equation2}$|_{\mu=\mu_0}$ must be lower-stable and upper-unstable.
        Therefore, $x=x_{2}(t;\lambda)$ and $x=x_{3}(t;\lambda)$ retain their hyperbolicity and stability for $\lambda\in(\lambda_{*},\lambda_{0}]$, and still exist (possibly coinciding) at $\lambda=\lambda_{*}$. However, again using Lemma \ref{LHopen}, $\eta_{\lambda_*,\mu_0}\in\partial\mathcal L_H\subset\mathcal L_1$.
        It follows from statement (i) of Proposition \ref{multiplicity of non-zero limit cycle} that the equation has at most one non-zero limit cycle, counted with multiplicity. This yields a contradiction. As a result, $m\leq2$ and the assertion follows.

        (i.2) We again consider the $\lambda$-parametric family \eqref{equation2}$|_{\mu=\mu_0}$ and follow the argument for contradiction as in statement (i.1).
        The only difference here is that the assumption $V_{2}(\mu_{0})<0$ implies that $x=x_{1}(t;\lambda_0)$ is unstable and $x=x_{2}(t;\lambda_0)$ is stable. Thus, based on Proposition \ref{property of rotated} and statement (ii) of Proposition \ref{multiplicity of non-zero limit cycle}, the limit cycles $x=x_{1}(t;\lambda)$ and $x=x_{2}(t;\lambda)$ approach each other as $\lambda$ decreases from $\lambda_0$, and both retain the hyperbolicity and stability for $\lambda\in(\lambda_{*},\lambda_{0}]$. This results in two positive limit cycles (counted with multiplicities) of equation \eqref{equation2}$|_{\lambda=\lambda_*,\mu=\mu_0}$, which again contradicts the known maximum number of non-zero limit cycles for the equation. Accordingly, equation \eqref{equation2}$|_{\lambda=\lambda_0,\mu=\mu_0}$ has at most one positive limit cycle. The conclusion of the statement is valid.


        (ii) 
        By assumption and Lemma \ref{LHopen}, $\eta_{\lambda,\mu}\in\mathcal L_H$ for $(\lambda,\mu)$ in some neighborhood of $(\lambda_0,\mu_0)$.
        Similar to the argument in statement (i.1), assume for a contradiction that equation \eqref{equation2}$|_{\lambda=\lambda_0,\mu=\mu_0}$ (i.e. equation \eqref{equation1}) possesses two consecutive positive limit cycles (which may coincide). Using Proposition \ref{property of rotated} and statement (ii) of Proposition \ref{multiplicity of non-zero limit cycle} if necessary, we know that for sufficiently small $\varepsilon>0$, equation \eqref{equation2}$|_{\lambda=\lambda_0-\varepsilon,\mu=\mu_0}$ has two hyperbolic positive limit cycles, with $V_3(\lambda_0-\varepsilon)\neq0$. Now consider the $\mu$-parametric family \eqref{equation2}$|_{\lambda=\lambda_0-\varepsilon}$.
        Note that $V_2(\mu)$ is non-trivial linear in $\mu$ and $V_2(\mu_0)=0$. One can slightly change the value of $\mu$ from $\mu_{0}$ such that $V_{2}(\mu) V_{3}(\lambda_0-\varepsilon)<0$. This changes the stability of $x=0$ for the family and then induces a Hopf-like bifurcation, generating a new positive limit cycle. In other words, there exists an appropriate small $\sigma$ such that equation \eqref{equation2}$|_{\lambda=\lambda_0-\varepsilon,\mu=\mu_0-\sigma}$ has at least three positive limit cycles, contradicting statement (i).

    (iii) It is sufficient to show that equation \eqref{equation2}$|_{\lambda=\lambda_0,\mu=\mu_0}$ (i.e. equation \eqref{equation1}) cannot have a positive limit cycle. In fact, if such a limit cycle exists and is denoted by $x=x_{1}(t;\lambda_{0},\mu_{0})$, then one of the cases below occurs.

    \noindent\text{ Case 1}: $x=x_{1}(t;\lambda_{0},\mu_{0})$ is non-hyperbolic.
    By Proposition \ref{property of rotated} and Proposition \ref{multiplicity of non-zero limit cycle}, it splits into two hyperbolic positive limit cycles of equation \eqref{equation2}$|_{\mu=\mu_0}$ when $\lambda$ slightly decreases from $\lambda_{0}$. Moreover, for this equation we have $V_{3}(\lambda)\not=0$ since it is non-trivial linear in $\lambda$.

    \noindent\text{ Case 2}: $x=x_{1}(t;\lambda_0,\mu_0)$ is hyperbolic.
    We can slightly change $\lambda$ from $\lambda_0$ such that $V_{3}(\lambda)V_{4}(\lambda_0,\mu_0)<0$. Then equation \eqref{equation2}$|_{\mu=\mu_0}$ has a hyperbolic limit cycle near $x=x_{1}(t;\lambda_0,\mu_0)$ and a new positive one bifurcating from $x=0$.

    In any case, there exists an appropriately small $\varepsilon$ such that equation \eqref{equation2}$|_{\lambda=\lambda_0-\varepsilon,\mu=\mu_0}$ admits at least two positive limit cycles, with $V_3(\lambda_0-\varepsilon)\neq0$. Moreover, due to Lemma \ref{LHopen}, this $\varepsilon$ can be chosen so as to also ensure $\eta_{\lambda_0-\varepsilon,\mu_0}\in\mathcal L_H$. Together with $V_2(\mu_0)=0$ from assumption, this comes back to the case in statement (ii) and therefore yields a contradiction. Accordingly, equation \eqref{equation2}$|_{\lambda=\lambda_0,\mu=\mu_0}$ has no positive limit cycles.

 The proof is complete.
\end{proof}

\subsection{Proof of Theorem \ref{theorem0}}
So far, thanks to Lemma \ref{V2V3}, we have established that equation \eqref{equation1} (i.e., equation \eqref{eq0} considered in Theorem \ref{theorem0}) has at most $4$ non-zero limit cycles under hypothesis $\bm{(H)}$. We will now finally prove Theorem \ref{theorem0} by optimizing this upper bound.

\begin{proof}[Proof of Theorem \ref{theorem0}]
    As explained at the beginning of this section, we only need to verify the assertion for equation \eqref{equation1} under hypothesis $\bm{(H)}$.
    For the case where the second Lyapunov constant $V_{2}\leq 0$, Lemma \ref{V2V3} tells us that equation \eqref{equation1} has at most two non-zero limit cycles. Hence, it remains to address the case $V_{2}> 0$.

    To this end, we replace the parameters $\lambda_0$ and $\mu_0$ in equation \eqref{equation1} by the variable ones $\lambda$ and $\mu$, and consider the family $\dot{x}=S(t,x;\lambda,\mu)$ given in \eqref{equation2}. The associated parameter vector is $\eta_{\lambda,\mu}=(\lambda,\lambda_{1},\lambda_{2},\mu,\mu_{1},\mu_{2})$, and the second, third and fourth Lyapunov constants are denoted by $V_2(\mu)$, $V_3(\lambda)$ and $V_4(\lambda,\mu)$, respectively. It is known by \eqref{Lyapunov constants-2} that
    \begin{align}\label{eq4.22}
      V_2(\mu)=(\mu-\mu_0)I_0+V_2(\mu_0),\indent V_3(\lambda)=(\lambda-\lambda_0)I_0+V_3(\lambda_0).
    \end{align}
    Since $I_0=\int_{0}^{T}f_{0}(t)\,dt>0$ by assumption, $V_2(\mu)$ and $V_3(\lambda)$ have unique zeros $\mu^*=\mu_0-\frac{V_2(\mu_0)}{I_0}$ and $\lambda^*=\lambda_0-\frac{V_3(\lambda_0)}{I_0}$, respectively. Moreover, one has
    \begin{align}\label{eq4.2}
      \mu^*-\mu_0<0, \text{ and } \lambda^*-\lambda_0\leq0 \text{ (resp. $>0$) when } V_3(\lambda_0)\geq0 \text{ (resp. $<0$)}.
    \end{align}

    In what follows, we show by comparison that equation \eqref{equation2}$|_{\lambda=\lambda_0,\mu=\mu_0}$ (i.e., equation \eqref{equation1}) has no limit cycles in either region $x>0$ or $x<0$ when $V_3(\lambda_0)\geq0$. The analogous argument applies when $V_3(\lambda_0)<0$ and is therefore omitted.

    First, from \eqref{equation2} we get
    \begin{equation*}
        S(t,x;\lambda,\mu)-S(t,x;\lambda_0,\mu_0)=(\lambda-\lambda_0)f_{0}x^{3}+(\mu-\mu_0)f_{0}x^{2}.
    \end{equation*}
    Let $\Gamma=\gamma_1\cup\gamma_2\cup\{\eta_{\lambda^*,\mu^*}\}$ be a polyline in the parameter space, where
    \begin{align*}
      \gamma_1=\left\{\eta_{\lambda_0,\mu}\left|\mu\in(\mu^*,\mu_0]\right.\right\},\ \ \
      \gamma_2=\left\{\eta_{\lambda,\mu^*}\left|\lambda\in(\lambda^*,\lambda_0]\right.\right\}.
    \end{align*}
    Then, for each pair $(\lambda,\mu)$ with $\eta_{\lambda,\mu}\in\Gamma$, we get by \eqref{eq4.2} that
    \begin{align}\label{eq4.21}
      S(t,x;\lambda,\mu)\leq S(t,x;\lambda_0,\mu_0)\text{ for }
      \left\{
      \begin{aligned}
            &(t,x)\in[0,T]\times\mathbb R ,& \text{if }\eta_{\lambda,\mu}\in \gamma_1, \\
            &(t,x)\in[0,T]\times\mathbb R^+ ,& \text{if } \eta_{\lambda,\mu}\in \Gamma\setminus\gamma_1. \\
        \end{aligned}
      \right.
    \end{align}

    Next, note that $\eta_{\lambda_0,\mu_0}\in\Gamma\cap\mathcal L_H$. Since $\mathcal L_H$ is open from Lemma \ref{LHopen}, we may distinguish four subcases below.

    \noindent\text{ Subcase 1}: $\Gamma\subset\mathcal L_H$. In this situation, $\eta_{\lambda^*,\mu^*}\in\mathcal L_H$. By statement (ii) of Proposition \ref{Lyapunov constants} and statement (iii) of Lemma \ref{V2V3}, equation \eqref{equation2}$|_{\lambda=\lambda^*,\mu=\mu^*}$ has a unique limit cycle at $x=0$, with multiplicity four and $V_4(\mu^*)>0$. Taking into account \eqref{eq4.21}, any solution $x=x(t)$ of equation \eqref{equation2}$|_{\lambda=\lambda_0,\mu=\mu_0}$ with initial value $x(0)>0$, satisfies $x(T,x_0)>x_0$ when it is well-defined on $[0,T]$. Hence, equation \eqref{equation2}$|_{\lambda=\lambda_0,\mu=\mu_0}$ has no limit cycles in the region $x>0$.

    \noindent\text{ Subcase 2}: $\eta_{\lambda^*,\mu^*}\in\partial\mathcal L_H$. Due to Lemma \ref{LHopen}, $\eta_{\lambda^*,\mu^*}\in\mathcal L_1$. We know by statement (i) Proposition \ref{Lyapunov constants} that equation \eqref{equation2}$|_{\lambda=\lambda^*,\mu=\mu^*}$ has a center at $x=0$. Similar to the analysis in Subcase 1, by comparison equation \eqref{equation2}$|_{\lambda=\lambda_0,\mu=\mu_0}$ has no limit cycles in the region $x>0$.

    \noindent\text{ Subcase 3}: $\gamma_2\cap\partial\mathcal L_H\not=\emptyset$. There exists $\lambda'\in(\lambda^*,\lambda_0]$ such that $\eta_{\lambda',\mu^*}\in \gamma_2\cap\partial\mathcal L_H\subset\mathcal L_1$. Observe that $V_3(\lambda')> V_3(\lambda^*)=0$ by \eqref{eq4.22}. From statement (i) of Proposition \ref{Lyapunov constants}, $x=0$ is the unique limit cycle of equation \eqref{equation2}$|_{\lambda=\lambda',\mu=\mu^*}$, and it is unstable. Again using \eqref{eq4.21} and comparison, equation \eqref{equation2}$|_{\lambda=\lambda_0,\mu=\mu_0}$ has no limit cycles in the region $x>0$.

    \noindent\text{ Subcase 4}: $\gamma_1\cap\partial\mathcal L_H\not=\emptyset$. There exists $\mu'\in(\mu^*,\mu_0)$ such that $\eta_{\lambda_0,\mu'}\in \gamma_1\cap\partial\mathcal L_H\subset\mathcal L_1$.  According to statement (i) of Proposition \ref{multiplicity of non-zero limit cycle}, equation \eqref{equation2}$|_{\lambda=\lambda_0,\mu=\mu'}$ admits at most one non-zero limit cycle. Moreover, since one has $V_2(\mu')>V_2(\mu^*)=0$, $x=0$ is an upper-unstable and lower-stable limit cycle of equation \eqref{equation2}$|_{\lambda=\lambda_0,\mu=\mu'}$.
    Therefore, in either region $x>0$ or $x<0$, any solution $x=x(t)$ of the equation well-defined on $[0,T]$ satisfies $x(T)>x_0$. Combining this with \eqref{eq4.21}, equation \eqref{equation2}$|_{\lambda=\lambda_0,\mu=\mu_0}$ has no limit cycles in either region $x>0$ or $x<0$.

    Based on the above argument and statement (i.1) of Lemma \ref{V2V3}, we verify that the total number of non-zero limit cycles of equation \eqref{equation2}$|_{\lambda=\lambda_0,\mu=\mu_0}$ (i.e., equation \eqref{equation1}) is at most two.

    Here we briefly show that this upper bound can be achieved. Consider equation \eqref{equation2}$|_{\lambda=I_1,\mu=I_2}$ with a specific parameter vector $\eta_{I_1,I_2}=(I_1,-I_0,0,I_2,0,-I_0)$, i.e.,
    \begin{align*}
      \frac{dx}{dt}=(I_1f_0-I_0f_1)x^3+(I_2f_0-I_0f_2)x^2,
    \end{align*}
    where $I_i$'s are defined as in \eqref{Lyapunov constants-2'}.
    It is easy to check by \eqref{Lyapunov constants-2} that $V_2(I_2)=V_3(I_1)=0$. Moreover, as shown in the proof of statement (i) of Proposition \ref{Lyapunov constants}, once $\eta_{I_1,I_2}\in\mathcal L_1$, the coefficients $I_1f_0-I_0f_1$ and $I_2f_0-I_0f_2$ of the equation are linearly dependent. This contradicts the linear independence of $f_0$, $f_1$ and $f_2$. Hence, $\eta_{I_1,I_2}\not\in\mathcal L_1$, and it follows from statement (ii) of Proposition \ref{Lyapunov constants} that $V_4(I_1,I_2)\neq0$. That is, $x=0$ is a limit cycle with multiplicity four. Accordingly, one can slightly change the parameters $\lambda$ and $\mu$ from $I_1$ and $I_2$, respectively, and then generate two non-zero limit cycles of equation \eqref{equation2} from $x=0$ by Hopf-like bifurcations.

The proof is complete.
\end{proof}

\section{Proofs of Theorems \ref{app1}, \ref{main theorem1} and \ref{main theorem}}
In this section, we apply Theorem \ref{theorem0} to the Smale-Pugh problem for the Abel equations \eqref{equationa}, \eqref{main equation} and \eqref{main equation1}, and then prove Theorems \ref{app1}, \ref{main theorem1} and \ref{main theorem}, respectively.
The key in these applications is to determine whether each class of coefficients can be characterized by an ET-system $\{f_{0},f_{1},f_{2}\}$ on $[0,T)$, with $f_{0}\not=0$. To begin with, we introduce some additional notions and a result related to Chebyshev systems. For more details of this theory the reader is referred to \cite{KS,K}.

As mentioned in Section 1, a set of functions $\mathcal F=\{f_{0},f_{1},\cdots,f_{m}\}$ forms an ET-system on an interval $E$, if the functions in $\mathcal F$ belong to $C^{m}(E)$ and the maximum number of zeros (counted with multiplicities) of their non-trivial combinations is $m$.
Now we additionally say that the ordered set $\mathcal F$ forms an \textit{extended complete Chebyshev system} (ECT-system) on $E$, if for each $k\in \{0,1,\cdots,m\}$, the subset $\{f_{0},f_{1},\cdots,f_{k}\}$ is an ET-system on $E$.
Clearly, an ECT-system is naturally an ET-system.
Moreover, the following classical result given in \cite{KS} enables us to determine an ECT-system using continuous Wronskians.

\begin{definition}
    Let $f_{0}, \cdots,f_{k}$ be $C^{k}$ functions on an interval $E\subset \mathbb{R}$.
    The continuous Wronskian of $\{f_{0}, \cdots,f_{k}\}$ at $t\in E$ is
    \begin{equation*}
W_{f_{0}, \cdots,f_{k}}(t) = \det(f_j^{(i)}(t); \ 0 \leq i, j \leq k) =
\begin{vmatrix}
    f_0(t) & \cdots & f_k(t) \\
    f_0'(t) & \cdots & f_k'(t) \\
    \vdots & \ddots & \vdots \\
    f_0^{(k)}(t) & \cdots & f_k^{(k)}(t)
\end{vmatrix}.
    \end{equation*}
\end{definition}

\begin{lemma}(\cite{KS})\label{KS1}
    Let \( f_0, f_1, \dots, f_m \) be $C^{m}(E)$ functions on an interval $E\subset \mathbb R$.
 The ordered set \( \{ f_0, f_1, \dots, f_m \} \) is an ECT-system on \( E \) if and only if, for each \( k = 0, 1, \dots, m \),
    \[
    W_{f_{0}, \cdots,f_{k}}(t) \neq 0 \quad \text{for all } t \in E.
    \]
\end{lemma}




We are ready to prove Theorem \ref{app1}, Theorem \ref{main theorem1} and Theorem \ref{main theorem}.
\begin{proof}[Proof of Theorem \ref{app1}]
Applying Theorem \ref{theorem0}, the assertion is true once $\{1,\sin t,\cos t\}$ forms an ET-system on $[0,2\pi)$.
This is clear because for any $(\lambda_{1},\lambda_{2},\lambda_{3})\neq(0,0,0)$,
\begin{align*}
\lambda_{1}+\lambda_{2}\sin t+\lambda_{3}\cos t=
  \left\{
      \begin{aligned}
            &\lambda_{1}+\sqrt{\lambda^{2}_{2}+\lambda^{2}_{3}}\sin \left(t+\arctan\frac{\lambda_{3}}{\lambda_{2}}\right) ,& \text{if }\lambda_2\neq0, \\
            &\lambda_{1}+\lambda_{3}\cos t ,& \text{if } \lambda_2=0, \\
        \end{aligned}
      \right.
\end{align*}
which has at most two zeros (counted with multiplicity) on $[0,2\pi)$. In addition, this upper bound is sharp when $|\lambda_{1}|\leq \sqrt{\lambda^{2}_{2}+\lambda^{2}_{3}}$.
\end{proof}

\begin{proof}[Proof of Theorem \ref{main theorem1}]
    By Theorem \ref{theorem0}, it is sufficient to show that $\{1,t,t^{2}\}$ forms an ET-system on $[0,1)$. This is trivial and also follows from Lemma \ref{KS1} and the equality $2W_{1}(t)=2W_{1,t}(t)=W_{1,t,t^{2}}(t)=2$.
\end{proof}

    \begin{proof}[Proof of Theorem \ref{main theorem}]
    Due to the time-rescaling $t\mapsto t^{\frac{1}{m_0+1}}$ and the arbitrariness of the parameters $a_{i}$ and $b_{i}$ with $i=0,1,2$,
    the assertion for equation \eqref{main equation1} can be verified by studying the following equation
    \begin{equation*}\label{40}
        \frac{dx}{dt}=\left(a_{0}+a_{1}t^{\alpha}+a_{2}t^{\beta}\right)x^{3}+\left(b_{0}+b_{1}t^{\alpha}+b_{2}t^{\beta}\right)x^{2},\indent t\in [0,1],
    \end{equation*}
    where $\alpha=\frac{m_{1}-m_{0}}{m_{0}+1}>0$ and $\beta=\frac{m_{2}-m_{0}}{m_{0}+1}>0$.
    Furthermore, taking the change of variable $t\rightarrow 1-t$ and absorbing the signs into the parameters, the above equation becomes
        \begin{equation}\label{41}
        \begin{split}
        \frac{dx}{dt}=&\left(a_{0}+a_{1}(1-t)^{\alpha}+a_{2}(1-t)^{\beta}\right)x^{3}\\&+
        \left(b_{0}+b_{1}(1-t)^{\alpha}+b_{2}(1-t)^{\beta}\right)x^{2},\indent t\in [0,1].
        \end{split}
    \end{equation}

    It is clear that the coefficients of equation \eqref{41} are generated by $1$, $(1-t)^{\alpha}$ and $(1-t)^{\beta}$.
In addition, for $t\in[0,1)$ a direct calculation yields
    \begin{align*}
        \begin{split}
            &W_{1}(t)=1,\\
            &W_{1,(1-t)^{\alpha}}(t){=-\alpha (1-t)^{\alpha-1}<0},\\
            &W_{1,(1-t)^{\alpha},(1-t)^{\beta}}(t)=\alpha\beta(1-t)^{\alpha+\beta-3}>0.
        \end{split}
    \end{align*}
    Thus, we know by Lemma \ref{KS1} that the set $\{1,(1-t)^{\alpha},(1-t)^{\beta}\}$ forms an ECT-system on $[0,1)$, and therefore is naturally an ET-system on $[0,1)$.
    According to Theorem \ref{theorem0}, equation \eqref{41} has at most three limit cycles (including $x=0$), and this upper bound is sharp. The assertion follows.
\end{proof}

\section{{Appendix}}
{Following the same argument in \cite{YHL}, we present here the details of obtaining Proposition \ref{YHLproposition}.}

\begin{proof}[{Proof of Proposition \ref{YHLproposition}}]
    {For the sake of brevity, we write $x(t)$ for $x(t,x_0)$ with $x(0)=x_0$.
    According to \eqref{YHLeqs}, we have
    \begin{equation}\label{p1}
        P'(x_{0})=\exp\left(\int_{0}^{T}3A(t)x^{2}(t)+2B(t)x(t)\,dt\right)=1
    \end{equation}
and
    \begin{equation}\label{p2}
        P''(x_{0})=\int_{0}^{T}\left(6A(t)x(t)+2B(t)\right)\exp\left(\int_{0}^{t}\left(3A(s)x^{2}(s)+2B(s)x(s)\right)\,ds\right)\,dt.
    \end{equation}
Since $x(t)$ is a non-zero periodic solution, we get that
\begin{equation}\label{p3}
    \int_{0}^{T}A(t)x^{2}(t)+B(t)x(t)\,dt=\int_{0}^{T}\frac{dx(t)}{x(t)}=\ln\frac{x(T)}{x_{0}}=0.
\end{equation}
Then, equalities \eqref{p1} and \eqref{p3} yield
\begin{equation*}
    h(T)=\int_{0}^{T}A(t)x^{2}(t)\,dt=0.
\end{equation*}
This verifies the validity of equation \eqref{77}.}

{Next, we simplify the expression of $P''(x_{0})$ in \eqref{p2}.
Note that $P''(x_{0})$ can be decomposed into $P''(x_{0})=4W_{1}+2W_{2}$, where
\begin{equation*}
    \begin{split}
            W_1 &= \int_{0}^{T} A(t)x(t)\exp\!\left(\int_{0}^{t} \bigl(2B(s)x(s) + 3A(s)x^2(s)\bigr)\, ds \right) dt,\\
            W_2 &= \int_{0}^{T} \bigl(B(t) + A(t)x(t)\bigr)\exp\!\left(\int_{0}^{t} \bigl(2B(s)x(s) + 3A(s)x^2(s)\bigr)\, ds \right) dt.
    \end{split}
\end{equation*}
Furthermore, note that we have already known $h(T)=h(0)=0$. Therefore,
\begin{equation*}
    \begin{split}
    W_1 &= \int_{0}^{T} A(t)x(t)\exp\!\left( 2\ln\frac{x(t)}{x_{0}} + h(t) \right) dt= \frac{1}{x_{0}^{2}} \int_{0}^{T} A(t)x^{3}(t)\exp(h(t))\, dt\\
    &= \frac{1}{x_{0}^{2}} \int_{0}^{T} x(t)\, d(\exp h(t))= -\frac{1}{x_{0}^{2}} \int_{0}^{T} \exp(h(t))\, d x(t),\\
    W_{2}&= \int_{0}^{T} \left( B(t) + A(t)x(t) \right)\exp\!\left( 2\ln \frac{x(t)}{x_{0}} + h(t) \right) dt \\
    &= \frac{1}{x_{0}^{2}}\int_{0}^{T} \left( B(t)x^{2}(t) + A(t)x^{3}(t) \right)
    \exp(h(t)) \, dt \\
    &= \frac{1}{x_{0}^{2}} \int_{0}^{T} \exp(h(t)) \, d x(t).
    \end{split}
\end{equation*}
Accordingly, the expression of $P''(x_{0})$ can be reduced to
\begin{equation*}
            P''(x_{0})=-\frac{2}{x_{0}^{2}}\int_{0}^{T}\exp h(t)\,dx(t).
\end{equation*}
Equality \eqref{88} is obtained and the proof is finished.}
\end{proof}

\section*{Data availability}
No data was used for the research described in the article.

\section*{Acknowledgements}
The authors thanks the referees, professor Jos\'{e} Luis Bravo Trinidad and professor Changjian Liu  for useful comments and suggestions which help us to improve both mathematics and presentations of this paper.

The first author is supported by NNSF of China (No. 12271212). The second and third authors are supported by the NNSF of China (No. 12371183 and No. 124B2006).

\end{document}